\documentclass[10pt]{article}
\usepackage{amssymb}
\usepackage{amsmath}
\usepackage{amsthm}
\theoremstyle{plain}
\newtheorem{theorem}{Theorem} \newtheorem{lemma}{Lemma}[section]
\newtheorem{propo}{Proposition}[section]
\newtheorem{question}{Question}\newtheorem{example}{Example}
\newtheorem{corol}{Corollary}[section] \newtheorem{defin}{Definition}
[section]
 
\newtheorem{conjecture}{Conjecture} 
\newcommand{\e}{\epsilon}
\newcommand{\Z}{\mathbb{Z}}\newcommand{\R}{\mathbb{R}}
\newcommand{\N}{\mathbb{N}}\newcommand{\al}{\alpha}\newcommand{\C}{\mathbb{C}}
\newcommand{\cM}{\mathcal{M}}

\newcommand{\cA}{\mathcal{A}}
\newcommand{\mat}{\mbox{Mat}\,}
\newcommand{\rmod}{\mbox{R-Mod}}  \newcommand{\fmod}{\mbox{$KF_r$-Mod}}  
\newcommand{\qmod}{\mbox{$Q_r$-Mod}} 
\newcommand{\flmod}{\mbox{$KF_{r-1}$-Mod}}
\newcommand{\mods}{\mbox{Mod-S}} 
\newcommand{\rank}{\mbox{rank}} \newcommand{\rk}{rk} 
\newcommand{\Endo}{\mbox{End}}\newcommand{\Syl}{\mbox{Syl}}
\newcommand{\Fin}{\mbox{Fin}}\newcommand{\Rank}{\mbox{Rank}}
\newcommand{\Hom}{\mbox{Hom}}\newcommand{\Ker}{\mbox{Ker}\,}
\newcommand{\Tor}{\mbox{Tor}}\newcommand{\Id}{\mbox{Id}}
\newcommand{\Ima}{\mbox{Im}\,}
\newcommand{\Ind}{\mbox{Ind}}\newcommand{\Conv}{\mbox{Conv}}
\newcommand{\kmod}{\mbox{K[X]-Mod}}\newcommand{\Span}{\mbox{Span}}
\newcommand{\KR} {KF_r}
\newcommand{\pbe}{[p_1 N^i_n]}\newcommand{\pbk}{[p_2 N^i_n]}
\title{Infinite dimensional representations of finite dimensional algebras
and amenability.
\footnote{
\, Research partly sponsored by MTA Renyi ``Lendulet'' Groups
and Graphs Research Group}}
\author{G\'abor Elek}

\begin{document}
\maketitle
\begin{abstract}
We present a novel quantitative  approach to the representation theory of finite
dimensional algebras motivated by the emerging theory of graph limits.

We introduce the rank spectrum of a finite dimensional algebra $R$ over a 
finite field.
The elements of the rank spectrum are representations of the algebra 
into von Neumann regular rank algebras, and two representations are 
considered to be equivalent 
if they induce
the same Sylvester rank functions on $R$-matrices.
 
Based on this approach, we can divide the finite dimensional algebras into
three types: finite, amenable and non-amenable representation types.
We prove that string algebras are of amenable representation type, but 
the wild Kronecker algebras are not. Here, the amenability of
the rank algebras associated to the limit points in
the rank spectrum plays a very important part. 

We also show that the limit points
of finite dimensional representations of algebras of amenable 
representation type can always be viewed as representations of the 
algebra in the continuous
ring invented by John von Neumann in the 1930's. 

As an application in algorithm theory, we introduce and study the notion of parameter testing
of modules over finite dimensional
algebras, that is analogous to the testing of bounded 
degree graphs introduced by Goldreich and Ron. We shall see that for string
algebras all the reasonable (stable) parameters are testable.
\end{abstract}\vskip 0.2in
\noindent
\newpage
\tableofcontents
\newpage
\section{Introduction}
\noindent
By Maschke's theorem, the complex group algebra $\C G$ is semisimple if
$G$ is a finite group. Hence all the infinite dimensional representations
of $G$ can be decomposed into finite dimensional subrepresentations.
The situation becomes much more interesting if the characteristic
of the coefficient field divides the order of the group.
\begin{example}\label{example1} Let $G=\Z_2\times \Z_2$ be the Klein group 
and $K$ be
the field of two elements. Let $K(t)$ be the transcendent extension of
$K$ by the element $t$. Then
$$\phi(a)=\begin{pmatrix} 1 & 1 \\ 0 & 1 \end{pmatrix},\quad
\phi(b)=\begin{pmatrix} 1 & t \\ 0 & 1 \end{pmatrix}$$
\noindent
defines an indecomposable representation $\phi:G\to GL(2,K(t))$, where
$G=\{1,a\}\times\{1,b\}.$ That is, the representation of $G$ on the 
infinite dimensional
space $K(t)\times K(t)$ cannot be decomposed into finite (or even infinite)
dimensional subrepresentations. \end{example}

\noindent
Infinite dimensional representations
of finite dimensional algebras has been studied for decades
(see e.g. \cite{Ringel79}, \cite{Craw}). The aim of our paper is to 
develop a theory for such representations guided by the convergence
and limit theory of finite graphs (see the monograph of L\'aszl\'o Lov\'asz
\cite{Lovasz}). Our philosophy is to view infinite dimensional representations
of a given algebra (when it is possible) as a sort of limit of its
finite dimensional representations. 

\vskip 0.05in
\noindent
{\bf The rank spectrum.}
Let $K$ be a finite field and $R$ be a finite dimensional $K$-algebra.
Let $\rmod$ denote the set of finitely generated (left) $R$-modules 
(up to isomorphism). 
If $M\in\rmod$, then $\phi^M:R\to\Endo_K(M)\cong \mat_{\dim(M)\times\dim(M)}(K)$
is the corresponding representation, where $\phi^M(r)=rm$ and $\dim(M)$ is
the $K$-dimension of the module $M$.
 For any $k,l\geq 1$ the map $\phi^M$ extends
naturally to the homomorphism
$$\phi^M_{k,l}:\mat_{k\times l}(R)\to \mat_{k\dim(M)\times \,l\dim(M)}(K)\,.$$
\noindent If $A\in \mat_{k\times l}(R)$, then let 
$$\rk_M(A):=\frac{\rank(\phi^M_{k,l}(A))}{\dim(M)}\,.$$
\noindent
Then $\rk_M$ is a Sylvester rank function, that is
\begin{itemize}
\item $\rk_M(I)=1$, where $I$ is the unit element of $R$ considered as 
a $1\times 1$-matrix. 
\item $
\rk_M\begin{pmatrix} A & 0 \\ 0 & B \end{pmatrix}= \rk_M(A)+\rk_M(B)\,,$
\item
$\rk_M\begin{pmatrix} A & 0 \\ C & B \end{pmatrix}\geq  \rk_M(A)+\rk_M(B)\,,$
\item
$\rk_M(AB)\leq \min(\rk_M(A),\rk_M(B))$\,,
\end{itemize}
provided that the $R$-matrices $A,B,C$ have appropriate sizes.
Observe that
$$\rk_M=\rk_{M\oplus M\oplus\dots\oplus M}\,.$$
We will show that if for two modules $\rk_M=\rk_N$, 
then there exist $k,l\geq 1$ such that
$M^k\cong N^l$. We say that a sequence $\{M_n\}^\infty_{n=1}\subseteq \rmod$ is 
{\bf convergent} if for
all pairs $k,l\geq 1$ and $A\in \mat_{k\times l}(R)$, 
$\lim_{n\to\infty} \rk_{M_n}(A)$ exists. 
The notion of convergence is based
on the Benjamini-Schramm convergence notion 
(see the Introduction of \cite{Eleklinear}).
Since the set of all finite matrices over $R$ is countable, 
from each sequence $\{M_n\}^\infty_{n=1}$
one can pick a convergent subsequence. Let us consider the locally convex, 
topological vector space
$C=\R^{\mat(R)}$, where $\mat(R)$ is the countable set of all finite 
matrices over $R$. The set of all
Sylvester rank functions on $\mat(R)$, $\Syl(R)$ forms a convex, 
compact subset of $C$. We will also consider
the subset of Sylvester rank functions corresponding to 
finite dimensional representations
$$\Fin(R):=\{\rk_M\,\mid\, M\in\rmod\}\,.$$
\noindent
Let $K$ be a finite field as above. A $K$-algebra $S$ is called a 
{\bf rank algebra} (see \cite{Good} for an introduction of rank functions on
von Neumann regular rings) if
\begin{itemize}
\item $S$ is a von Neumann regular ring,
\item $\mat(S)$ is equipped with a Sylvester rank function $\rank_S$.
\item $\rank_S(A)=0$ implies $A=0$.
\end{itemize}
\noindent
The most important examples of rank algebras are skew fields, matrix rings
over skew fields and more generally, semisimple Artinian rings.
If we have a
representation $\rho:R\to S$ into a rank algebra we still have an 
associated Sylvester rank function
$$\rk_\rho(A):=\rank_S(\tilde{\rho}(A))\,,$$
where $\tilde{\rho}$ is the appropriate extension for matrices. We say that
two such representations $\rho_1:R\to S_1$ and $\rho_2:R\to S_2$ are
equivalent if $\rk_{\rho_1}=\rk_{\rho_2}$. 
The {\bf rank spectrum} of $R$, $\Rank(R)$ is the set of Sylvester rank
functions that are in the form $\rk_\rho$ for some representation $\rho:R\to S$.
We say that the infinite dimensional representation $\rho:R\to S$ is the {\bf limit}
of the convergent sequence of finite dimensional representations $\{\rho_n:R\to \mat_{m_n\times m_n}(K)\}^\infty_{n=1}$
if for any $A\in\mat(R)$, $$\lim_{n\to\infty} \rk_{\rho_n}(A)=\rk_\rho(A)\,$$
that is, if $\rk_\rho$ is the limit point of the sequence $\{\rk_{\rho_n}\}^\infty_{n=1}$ in $\Syl(R)$.
\begin{example} \label{example2}
Let $\{K(t_n)=E_n\}^\infty_{n=1}$ be a sequence of finite dimensional field extensions of $K$ with generators
$\{t_n\}^\infty_{n=1}$. Suppose that $\dim_K(E_n)\to\infty$. Let 
$G=\Z_2\times \Z_2$ as in Example \ref{example1} and define the representation $\phi_n:KG\to \mat_{2\times 2}(E_n)$ by
$$\phi_n(a)=\begin{pmatrix} 1 & 1 \\ 0 & 1 \end{pmatrix},\quad
\phi_n(b)=\begin{pmatrix} 1 & t_n \\ 0 & 1 \end{pmatrix}\,.$$
\noindent
Then, the representations $\{\phi_n\}^\infty_{n=1}$ converge to the representation
$\phi:KG\to \mat_{2\times 2}(K(t))$ described in Example \ref{example1}.
\end{example}

\begin{propo} \label{P1} For any finite dimensional algebra $R$,
the set $\Rank(R)$ is a convex, compact subspace of
$\Syl(R)$.
\end{propo}
\noindent
The following question is the somewhat analogous to the famous Connes
Embedding Problem \cite{Capr}.
\begin{question} \label{connes} (see {\bf Theorem \ref{P2}})
Is the closure of $\Fin(R)$ in $\Syl(R)$ equals to $\Rank(R)$? 
\end{question}
\noindent
Note that if $\rk\in \Syl(R)$ is a Sylvester rank function such
that the values of $\rk$ are in the set $\Z/d$ for some $d\in\N$, then
by Schofield's Theorem (Theorem 7.12 \cite{Scho}), there exists a skew field $D$ and
a representation $\rho:R\to \mat_{d\times d}(D)$ such that $\rk=\rk_{\rho}$
\vskip 0.2in
\noindent
{\bf Hyperfiniteness.} 
Let $R$ be a finite dimensional algebra as above and 
$\{M_n\}^\infty_{n=1}\subset \rmod$ be
a set of modules. We say that $\{M_n\}^\infty_{n=1}$ is {\bf hyperfinite},
if for any $\e>0$, there exists a finite family 
$\{N_1,N_2,\dots, N_t\}\subset \rmod$ so that
for any $n\geq 1$, there exists a submodule $P_n\subset M_n,
\dim_K(P_n)\geq (1-\e)\dim_K(M_n)$ and non-negative integers $\alpha_1,\alpha_2,
\dots,\alpha_t$ such that $P_n\cong \bigoplus_{i=1}^t N_i^{\alpha_i}$.
An algebra $R$ of infinite representation type (that is the set of
indecomposable modules $\Ind(R)$ in $\rmod$ is infinite) is called of
{\bf amenable representation type}, if $\rmod$ itself is a hyperfinite family.
In the course of the paper, we shall see (Proposition \ref{stringhyper})  that string algebras are
of amenable representation type and the wild Kronecker algebras ({\bf Theorem \ref{kronecker}}) are not.
These results suggest a relation between the tame/wild and the amenable/non-amenable
 dichotomies.

\vskip 0.1in
\noindent
Let $d_R$ be a metric on $\Syl(R)$ defining the compact topology.
Two modules $M,N\in\rmod$  are $\e$-close to each other
if there exist submodules $P\subset M, Q\subset N, P\cong Q$ such that
$\dim_K(P)\geq (1-\e)\dim_K(M)$, $\dim_K(Q)\geq (1-\e)\dim_K(N)$.
Motivated by graph theoretical results (Theorem 5. \cite {Elekhyper}),
we make the following conjecture.
\begin{conjecture} \label{pip} Let $R$ be an finite dimensional algebra of 
amenable 
representation type.
Then, for any $\e>0$ we have a $\delta>0$ such that for two elements 
$M,N\subset \rmod$, $1-\delta\leq \frac{\dim_K(M)}{\dim_K(N)}\leq 1+\delta$, 
$d_R(\rk_M,\rk_N)<\delta$, then
$M$ and $N$ are $\e$-close to each other.
\end{conjecture}
\noindent
The meaning of the conjecture is that if we have a finite dimensional
algebra of amenable representation 
type, then the almost isomorphism of two $R$-modules
of the same dimension can be decided by checking the ranks of finitely
many matrices. This observation leads to the 
representational theoretical analogue of the Constant-Time Property Testing
Algorithms developed by Goldreich and Ron for bounded degree graphs.
(see Section \ref{property}). The main result of our paper is the following.
\begin{theorem} \label{stringtheorem}
Conjecture \ref{pip} holds if $R$ is a string algebra.
\end{theorem}
\noindent
\vskip 0.1in
\noindent
{\bf Amenability of infinite dimensional algebras.}
Amenability and hyperfiniteness are intimately related notions in group
theory. The notion of amenability for algebras were introduced by Misha 
Gromov \cite{Gromov} and plays an important part in the study of the 
rank spectrum. 
\begin{defin}\label{def1.1}
An element $\rk_\rho$ in the rank spectrum of $R$ is amenable, if there
exists a hyperfinite sequence $\{M_n\}^\infty_{n=1}\subset \rmod$ such that
$\lim_{n\to\infty}\rk_{M_n}=\rk_{\rho}$.
\end{defin}
\noindent
The following conjecture is motivated by the theorem of Oded Schramm on hyperfinite graph sequences \cite{Oded}.
\begin{conjecture}\label{con2}
If $\rho:R\to S$ defines $\rk_\rho$, where $S$ is
an amenable rank algebra, then $\rk_\rho$ is hyperfinite. 
\end{conjecture}
\noindent
We confirm
this conjecture in the case, where $S=\mat_{l\times l}(D)$ and $D$ is an
amenable skew field ({\bf Theorem \ref{fontostetel}}).
For the converse, one can hope that the limits of hyperfinite sequences
can be represented by homomorphisms $\rho:R\to S$, where $S$ is
an amenable rank algebra. 

\vskip 0.1in
\noindent
{\bf The continuous algebra of John von Neumann.} It is well-known that any countable amenable group algebra has
a trace-preserving embedding into the unique separable hyperfinite $II_1$-factor. In the world of ranks the
role of the hyperfinite $II_1$-factor will be played by an algebra constructed by John von Neumann in the 1930's
\cite{vonNeumann}. The construction goes as follows.
We consider the following sequence of diagonal
embeddings.
$$K \to \mat_{2\times 2} (K) \to  \mat_{4\times 4} (K)
 \to \mat_{8\times 8}  (K)\to
\dots\,,$$
where $K$ is a finite field.
Then all the embeddings are preserving the normalized ranks. 
Hence the direct limit $\varinjlim \mat_{2^k\times 2^k} (\C)$ is a rank algebra.
The addition, multiplication
and the rank function extends to the metric completion $M_K$ of
the direct limit algebra. The resulting algebra $M_K$ 
is a simple continuous rank algebra (see also \cite{Good}). 
\begin{conjecture}
All the countable dimensional amenable rank algebras embeds into $M_K$. Conversely, all the countable
dimensional rank subalgebra of $M_K$ is amenable.
\end{conjecture}
\noindent
\begin{theorem} \label{embedding}
Let $\cM=\{M_n\}^\infty_{n=1}\subset \rmod$ be a hyperfinite convergent sequence and $\rk_\cM=\lim_{n\to\infty} \rk_{M_n}$.
Then there exists a representation $\rho:R\to M_K$ such that $\rk_\rho=\rk_\cM$.
\end{theorem}

\vskip 0.2in
\noindent
{\bf Acknowledgement.} The author is grateful for the hospitality of the
Bernoulli Center at the \'Ecole Polytechnique F\'ed\'erale Lausanne , 
where some of this work was carried out. 

\section{The continuous ultraproduct of rank algebras}
\label{conultra}
First, let us prove Proposition \ref{P1}.
Let $\{S_n\}^\infty_{n=1}$ be rank algebras and $\{\rho_n:R\to S_n\}^\infty_{n=1}$ be a convergent sequence in the rank 
spectrum $\Rank(R)$, that is, for any matrix $A\in\mat(R)$, 
$\lim_{n\to\infty}\rk_{\rho_n}(A)$ exists.
Let $\omega$ be a non-principal ultrafilter on the
natural numbers and $\lim_\omega$ be the associated ultralimit.
We define the {\bf continuous ultraproduct} $\prod_\omega S_n$, a
quotient of the classical ultraproduct, the following way.
Let $I\subset \prod^\infty_{n=1} S_n$ be the set of elements
$\{a_n\}^\infty_{n=1}$ such that 
$$\lim_\omega \rank_{S_n}(a_n)=0\,.$$
\noindent
Clearly, $I$ is a two-sided ideal in $\prod^\infty_{n=1} S_n$. Then
$\Sigma=\prod_\omega S_n:= \prod^\infty_{n=1} S_n/I$ is a 
von Neumann regular ring. Indeed, if
$[\{a_n\}^\infty_{n=1}]\in \prod_\omega S_n$ and $x_n\in S_n, a_nx_na_n=a_n$, then
$$[\{a_n\}^\infty_{n=1}][\{x_n\}^\infty_{n=1}][\{a_n\}^\infty_{n=1}]=
[\{a_n\}^\infty_{n=1}]\,.$$
That is each element of $\prod_\omega S_n$ has a pseudoinverse.
Now, if $A\in\mat_{k\times l}(\Sigma)=[\{A_n\}^\infty_{n=1}]$, then
let $\rank_{\Sigma}(A)=\lim_\omega \rank_{S_n}(A_n)\,.$  Then
$\rank_\Sigma$ is a Sylvester rank on $\Sigma$ and $\rank_\Sigma(a)\neq 0$ if
$a\in \Sigma$\,.

\noindent
The representation $\rho_\Sigma:R\to\Sigma$ is defined
by $\rho_\Sigma(r):=[\{\rho_n(r)\}^\infty_{n=1}]\,.$
Then, $\rk_\Sigma(A)=\lim_{n\to\infty} \rk_{\rho_n}(A)$ holds for
any matrix $A\in\mat(R)$.
Hence, $\Rank(R)$ is a closed subset of $\Syl(R)$. 

\noindent
If the elements $\rk_1$ resp. $\rk_2$ are represented by homomorphisms
$\rho_1:R\to S$ and $\rho_2:R\to T$, where $S,T$ are rank algebras, then
$\lambda\cdot\rk_1+\mu\cdot\rk_2$,\\($\lambda,\mu\geq 0, \lambda+\mu=1$), is 
represented by $\rho:R\to S\oplus T$,
where
$$\rank_{S\oplus T}=\lambda\cdot\rank_S+\mu\cdot\rank_T\,.$$
\noindent
Hence, $\Rank(R)$ is a convex subset of $\Syl(R)$. \qed \vskip 0.2in
\noindent
The following lemma is implicite in \cite{Eleklinear}, nevertheless we include 
the proof for the sake of completeness.
\begin{lemma}
Let $K$ be a finite field, $R$ be a finite dimensional $K$-algebra and
$\rho:R\to S$ be a representation of $R$ into a rank algebra.
Then there exists a countable dimensional subalgebra $S'\subset S$ that
contains the image of $\rho$. That is, all the elements of the rank spectrum
can be represented by countable dimensional rank algebras.
\end{lemma}
\proof
Let $X$ be a finite dimensional linear subspace of $S$ containing
the unit. Denote by $P(X)$ the finite dimensional linear subspace
of $S$ spanned by elements in the form
$$\{x_1x_2\,\mid\, x_1,x_2\in X\}\,.$$
Also, let $R(X)$ be an arbitrary finite set such that
for any $x\in X$, there exists $y\in R(X)$ such that $xyx=x$. Let
$X_1=\mbox{Im}\,\rho(R)$ and $X_{n+1}=R(P(X_n))$. Then $\cup^\infty_{n=1} X_n$
is a countable dimensional von Neumann regular ring containing the image
of $R$. \qed
\section{The Ziegler spectrum}
The goal of this section is to recall the notion of the Ziegler spectrum
\cite{Zieg} (see also \cite{Prest}) and show that finitely generated 
indecomposable modules
are isolated points of the closure of $\Fin(R)$.

\noindent
Let $R$ be a finite dimensional algebra over a finite field $K$. A pp-formula
$\phi$ of type $t$ is given by a matrix 
$$A=\{a_{ij}\}_{1\leq i\leq m,1\leq j \leq n}\in \mat_{m\times n}(R)\,.$$
\noindent
Let $M$ be a module over $R$. We say that $\underline{v}=\{v_1,v_2,\dots,v_t\}
\in M^t$ satisfies $\phi$ if there exists $\{y_{t+1},y_{t+2},\dots,y_n\}\in M^{n-t}$
such that for any $1\leq i \leq m$
$$\sum_{j=1}^t a_{ij}v_j +\sum_{j=t+1}^n a_{ij} y_j=0\,.$$
The elements $\underline{v}$ satisfying $\phi$ form the $pp$-subspace
$M(\phi)\subset M^t\,.$
We say that two formulas (of the same type) are equivalent,
$\phi\cong \psi$, if for any $M$, $M(\phi)= M(\psi)$. Also,
$\phi\geq \psi$ if for any $M$, $M(\phi)\supset M(\psi)$.
For any $t\geq 1$, the equivalence classes of pp-formulas of type $t$
form the modular lattice $pp^t(R)$.
The Ziegler spectrum of $R$, $Zg(R)$, consists of the pure-injective
indecomposable modules over $R$. The basic open sets of $Zg(R)$
are given by pp-pairs $\langle\phi,\psi\rangle$, where
$\phi\geq \psi$.
$$U_{\phi,\psi}=\{M\,\mid\, M(\phi)\supsetneq  M(\psi)\}\,.$$

\noindent
Then, the set of isolated points in $Zg(R)$ are exactly the finitely
generated indecomposable modules and
they form a dense subset in the quasi-compact (not necessarily Hausdorff)
space $Zg(R)$ (Corollary 5.3.36. and Corollary 5.3.37. \cite{Prest}). Let $\phi$ be a $pp$-formula of type $t$ and
let $M$ be a finitely generated module.
Set
$$D_M([\phi]):=\frac{\dim_K(M(\phi))}{\dim_K(M)}\,.$$
\noindent
Then, $D_M$ is a dimension funcion on $pp^t(R)$, that is
\begin{itemize}
\item $D_M(0)=0$
\item $D_M(1)=1$
\item $D_M(a\wedge b)+ D_M(a \vee b)= D_M(a)+D_M(b)$
\item $D_M(a)\leq D_M(b)$ if $a\leq b$
\end{itemize}

\noindent
(Note that we do not require that $ a \lneq b => D_M(a) \lneq D_M(b)$\,.)
Now let $\rho:R\to S$ be a representation (possibly infinite dimensional) into
a rank algebra. Recall \cite{Good} that the finitely generated
right submodules of $S$, $\mods$, are projective, all exact sequences
of such modules split and the rank function defines a dimension function
$\dim_S$ on the modular lattice $\mods$. 
Let $A:S^m\to S^n$ be a module endomorphism (we nonchalantly consider $A$ as 
an $m\times n$-matrix in $\mat(S)$). Then 
$$\rank_S(A)=\dim_S(\Ima(A))=n-\dim_S(\Ker(A))\,.$$
\noindent
Let us consider $S$ as a left $R$-module.
\begin{lemma} \label{compute}
For any $\phi\in pp^t(R)$, $S(\phi)\in \mods$ and its
dimension $\dim_S(S(\phi))$ can be computed
by knowing only the element $\rk_\rho$ in the rank spectrum.
\end{lemma}
\proof
Recall that $S(\phi)$ is the set of elements
$\{s_1,s_2,\dots,s_t\}\in S^t$ such that there exists
$\{y_{t+1},y_{t+2},\dots,y_n\}\in S^{n-t}$ so that
for any $1\leq i \leq m$
$$\sum_{j=1}^t a_{ij}s_j+\sum_{j=t+1}^n a_{ij}y_j=0\,,$$
\noindent
where the matrix $A$ defines $\phi$. Therefore, $S(\phi)$ is a right
$S$-module. Let $Q$ be the set of elements 
$\{z_1,z_2,\dots,z_n\}\in S^n$ such that
$$\sum_{j=1}^n a_{ij}z_j=0\,,$$
\noindent
and $z_i=0$, if $1\leq i \leq t$.
Then
$$0\to Q\to \Ker(A)\to S(\phi)\to 0$$
\noindent
is an exact sequence, where $Q$ is the kernel of a matrix $B\in\mat(R)\subset
\mat(S)\,.$
That is $Q$ and $\Ker(A)$ are both finitely generated right $S$-modules
and $S(\phi)$ is also finitely generated since it is a quotient of
$\Ker(A)$. We have that
$$D_S(\phi)=\dim_S(S(\phi))=\dim_S(Q)-\dim_S(\Ker(A))=t+\rk_\rho (B)-
\rk_\rho (A)\quad\qed $$
\noindent
\begin{corol} \label{ppcorol}
If $\{M_n\}^\infty_{n=1}\in\rmod$ is a convergent sequence, then for any $pp$-formula $\phi$,
$\{D_{M_n}(\phi)\}^\infty_{n=1}$ converges. \end{corol}
So, one can identify the elements of the rank spectrum with certain
dimension functions on $pp^t(R)$. Note that if
$\{M_n\}^\infty_{n=1}$ is a convergence sequence of modules,
$\langle \phi,\psi \rangle$ is a pp-pair, for
any $n\geq 1$ $M_n(\phi)\neq M_n(\psi)$ but
$$\lim_{n\to\infty} (D_{M_n}(\phi)-D_{M_n}(\psi))=0\,.$$
\noindent 
then if $\rho:R\to S$ represents the limit of $\{M_n\}^\infty_{n=1}$
in the rank spectrum, we have $S(\phi)=S(\psi)$.
This phenomenon signifies an important difference between the
topology of the Ziegler-spectrum and the rank spectrum.
\section{The geometry of the rank spectrum}
The goal of this section is to understand the geometry of the closure
of $\Fin(R)$ (the ranks associated to finitely generated modules) in the
rank spectrum.
Let $M\in\rmod$. By the Krull-Schmidt theorem, $M$ can be written
in a unique way as $M\cong \oplus^s_{i=1} Q_i^{n_i}\,,$ where the $Q_i$'s are
indecomposable modules. The {\bf weight} of an indecomposable module
$Q_i$ in $M$ is defined as
$$w_{Q_i}(M):=\frac{n_i\dim_K(Q_i)}{\dim_K(M)}\,.$$
That is,
$$M\cong \bigoplus_{Q\in\Ind(R)} Q^{\frac{w_Q(M) \dim_K(M)}{\dim_K(Q)}}\,$$
\noindent
and
$\sum_{Q\in\Ind(R)} w_Q(M)=1\,.$
\begin{propo} \label{isomodules}
If $M,N\in\rmod$, then
$\rk_M=\rk_N$ if and only if $w_Q(M)=w_Q(N)$ for all $Q\in\Ind(R)$.
In particular, $\rk_Q\neq \rk_P$, whenever $P\neq Q\in \Ind(R)\,.$
\end{propo}
\proof
First, suppose that $M\cong \oplus^s_{i=1} Q_i^{n_i}$.
\begin{lemma}\label{sulyok}
$\rk_M=\sum^s_i w_{Q_i}(M) \rk_{Q_i}$
\end{lemma}
\proof
Let $A\in \mat_{m\times n}(R)$ be a matrix. Then, using the definition in
the Introduction,
$$\rank(\phi^{M}_{m,n}(A))=\sum^s_{i=1} w_{Q_i}(M) 
\frac{\dim_K(M)}{\dim_K(Q_i)} \rank(\phi^{Q_i}_{m,n}(A))\,.$$
\noindent
That is,
$$\rk_M(A)=\frac{\rank(\phi^{M}_{m,n}(A))}{\dim_K(M)}=
\sum^s_{i=1} w_{Q_i}(M)\rk_{Q_i}(A)\,.$$
\noindent
So, by the lemma above, if for all $Q\in\Ind(R), w_Q(M)=w_Q(N)$ holds, then
$\rk_M=\rk_N$. Note that it means that for certain integers $k$ and $l$,
$M^k\cong N^l$. \qed
\vskip 0.2in
\noindent
Now, suppose that for some $Q$, $w_Q(M)\neq w_Q(N)$. 
Let $\langle \phi,\psi \rangle$ be a pp-pair, that isolates $Q$, that is
\begin{itemize}
\item $Q(\phi)\neq Q(\psi)$
\item $P(\phi)=P(\psi)$ if $P\ncong Q, P\in \Ind(R)$\,.
\end{itemize}
\noindent
By, Lemma 1.2.3. \cite{Prest}, $M(\phi)=\bigoplus_{Q\in\Ind(R)} Q^{\frac{w_Q(M) \dim_K(M)}{\dim_K(Q)}}(\phi)\,.$
Hence,
\begin{equation}
\label{weight}D_M(\phi)-D_M(\psi)=w_Q(M)(D_Q(\phi)-D_Q(\psi))
\end{equation}
\noindent
That is,
$$D_M(\phi)-D_M(\psi)=w_Q(M)(D_Q(\phi)-D_Q(\psi))\neq
D_N(\phi)-D_N(\psi)\,.$$
\noindent
Hence by Lemma \ref{compute}, $\rk_M\neq \rk_N$.\,\qed
\vskip 0.2in
\noindent
Note that (\ref{weight}), has the following corollary.
\begin{corol} \label{weightcorollary}
If $\{M_m\}^\infty_{n=1}\subset \rmod$ is a convergent sequence, then
for all $Q\in\Ind(R)$ the sequence $\{w_{Q}(M_n)\}^\infty_{n=1}$ is
convergent as well. 
\end{corol}
\vskip 0.2in
\noindent
Now, let $Q_1,Q_2,\dots$ be an enumeration of the indecomposable
modules in $\rmod$.
Let $\{a_n\}^\infty_{n=1}$ be non-negative real numbers, such that
$\sum^\infty_{n=1} a_n=1$. Then $\underline{a}=
\sum^\infty_{n=1} a_n\rk_{Q_n}$
is a well-defined element of $\Rank(R)$. We denote these set of elements
by $\Conv(R)$.
\begin{propo}\label{convex}
Let $\{\underline{a^i}=\sum^\infty_{n=1} a^i_n \rk_{Q_n}\}^\infty_{i=1}$ be
a convergent sequence such that $\lim_{i\to\infty} \underline{a^i}=
\underline{b}=\sum^\infty_{n=1} b_n\rk_{Q_n}\in \Conv(R)$. Then for
any $n\geq 1$, $\lim_{i\to\infty} a^i_n=b_n$.
\end{propo}
\proof
Let suppose that for some $n\geq 1$ and subsequence $\{i_k\}^\infty_{k=1}$
$$\lim_{k\to\infty} a^{i_k}_n=c_n\neq b_n\,.$$
\noindent
Let $\langle \phi,\psi \rangle$ be a pp-pair that isolates $Q_n$.
Then,
$$\lim_{k\to\infty} (D_{\underline{a}^{i_k}}(\phi)-D_{\underline{a}^{i_k}}(\psi)
=c_n(D_{Q_n}(\phi)-D_{Q_n}(\psi))\neq $$
$$\neq b_n(D_{Q_n}(\phi)-D_{Q_n}(\psi))=D_{\underline{b}}(\phi)-
D_{\underline{b}}(\psi)\,,$$
leading to a contradiction. \qed
\begin{corol}
The limit points of $\Ind(R)$ in the rank spectrum form a non-empty
closed set that is disjoint from $\Conv(R)$.
\end{corol}
\noindent
\section{Sofic algebras}
Let $R$ be finite dimensional algebra over the finite field $K$ as in
the previous sections and $\rho:R\to S$ be a homomorphism to a rank algebra.
In this section we will get a sufficient condition under which $\rk_\rho$ is
an element of the closure of $\Fin(R)$ (see Question \ref{connes}).
Let $\cA$ be a countable dimensional algebra over the finite field $K$
with basis $\{1=e_1,e_2,\dots\}$. Following Arzhantseva and Paunescu \cite{Arzh},
we call $\cA$ a {\bf sofic algebra} 
\begin{itemize}
\item
if there exists a non-negative
function $j:\cA\to\R$ such that $j(a)>0$ if $a\neq 0$, and a sequence of real numbers $k_1> k_2>\dots$
tending to zero,
\item for any $n\geq 1$, there exists a unital, linear map
$\phi_n:\cA\to\mat_{m_n\times m_n}(K)$ so that
$$\frac{\Rank(\phi_n(ab)-\phi_n(a)\phi_n(b))}{m_n}<k_n\,,$$
whenever $a,b\in\mbox{Span}(e_1,e_2,\dots,e_n)$
\item $\frac{\Rank(\phi_n(a))}{m_n}>\frac{j(a)}{2}\,,$ if $n$ is large enough.
\end{itemize}
\noindent
We call such a sequence $\{\phi_n\}^\infty_{n=1}$ a sofic representation sequence of $\cA$.
It is not hard to see (Proposition 5.1 \cite{Eleklinear}) that $\cA$ is
sofic if and only if there exists an injective, unital homomorphism
$\phi:\cA\to \mat^{\cM}_\omega$, where $\mat^{\cM}_\omega$ is the continuous
ultraproduct of the matrix algebras $\cM:=\{\mat_{m_n\times m_n}(K)\}^\infty_{n=1}$.
Recall from Section \ref{conultra} that the rank $\rk_{\cM}$ on
$\mat^{\cM}_\omega$ is given by $\rk_{\cM}([\{a_n\}]=\lim_\omega \rk_{m_n}(a_n)$,
where $\rk_{m_n}$ is the normalized rank function on the matrix algebra
$\mat_{m_n\times m_n}(K)$.
By the proof of Proposition \ref{P1}, it follows that any element of
the closure of $\Fin(R)$ is associated to a unital homomorphism 
$\rho:R\to\mat^{\cM}_\omega$. Now we prove the converse.
\begin{theorem} \label{P2}
Let $\cA$ be a sofic algebra over our finite base field $K$ with a rank
$\rk_{\cA}$ derived from the injective homomorphism 
$\phi:\cA\to \mat^\cM_\omega$. Let $R$ be a finite dimensional algebra over
a finite field $K$, and $\rho:R\to\cA$ be the corresponding element of
the rank spectrum $\Rank(R)$. Then $\rk_\rho$ is in the closure of $\Fin(R)$.
\end{theorem}
\proof By Proposition 5.1 \cite{Eleklinear}, we have a sequence
of unital maps $\{\phi_n:\cA\to\mat_{m_n\times m_n}(K)\}^\infty_{n=1}$ such
that for any $\e>0$ and for any $a,b\in \cA$
$$\{n\,\mid\,\rk_{m_n}(\phi_n(ab)-\phi_n(a)\phi_n(b))<\e\}\in\omega$$
\noindent
and for any $k,l\geq 1$ and $A\in\mat_{k,l}(\cA)$
$$\{n\,\mid\,\frac{\dim_K(\phi^{k,l}_n(A))}{m_n}-\rk_\rho(A)|<\e\}\in\omega\,,$$
\noindent
where $\phi^{k,l}_n$ is the extension of $\phi_n$ onto $\mat_{k,l}(\cA)$.
Therefore, by taking a subsequence we can assure that for any $a,b\in \cA$
\begin{equation} \label{1nov15}
\lim_{n\to\infty}\rk_{m_n}(\phi_n(ab)-\phi_n(a)\phi_n(b))=0
\end{equation}
\noindent
and for any $k,l\geq 1$ and $A\in\mat_{k,l}(\cA)$
\begin{equation} \label{2nov15}
\lim_{n\to\infty}\frac{ \dim_K(\phi^{k,l}_n(A))}{m_n}=\rk_\rho(A).
\end{equation}
\noindent
If the sofic representation sequence $\{\phi_n\}^\infty_{n=1}$ satisfies    (\ref{2nov15}), then we call it a
{\bf convergent sofic representation sequence}.
It is enough to construct for all $n\geq 1$ a subspace $V_n\subset K^{m_n}$
such that
\begin{itemize}
\item $\lim_{n\to\infty}\frac{ \dim_K(V_n)}{m_n}=1\,.$
\item For any $a\in R$, $\phi'_n(a)(V_n)\subset V_n\,,$
\end{itemize}
\noindent
where $\phi_n'=\phi_n\circ \rho.$
Indeed, for such sequence of subspaces $\{V_n\}^\infty_{n=1}$ the maps
$\phi'_n$ define $R$-module structures on $V_n$ and for any matrix
$A\in\mat_{k,l}(\cA)$,

$$\lim_{n\to\infty}\left(\frac{\dim_K(\phi_n^{k,l}(A))}{m_n}-
\frac{\dim_K(\psi_n^{k,l}(A))}{m_n}\right)=0,$$
where $\psi_n^{k,l}$ is the restriction of $\phi_n^{k,l}$ onto $V^l_n\,.$
By (\ref{1nov15}), for any pair $a,b\in R$,
$$\lim_{n\to\infty} \frac{\dim_K(\Ker(\phi_n'(ab)-\phi_n'(a)\phi_n'(b))}{m_n}=0\,.$$
\noindent
Hence, for $V_n=\cap_{a,b\in R}\Ker(\phi_n'(ab)-\phi_n'(a)\phi_n'(b))$,
$\lim_{n\to\infty}\frac{ \dim_K(V_n)}{m_n}=1\,.$
We need to show that if $v\in V_n$, then $ \phi_n'(c)(v)\in V_n$, that is
for any $a,b\in R$
\begin{equation} \label{3nov15}
\phi_n'(ab)\phi_n'(c)(v)=\phi_n'(a)\phi_n'(b)\phi_n'(c) (v)\,.
\end{equation}
\noindent
Since $\phi_n'(ab)\phi_n'(c)(v)=\phi_n'(abc)(v)$ and
$\phi_n'(a)\phi_n'(b)\phi_n'(c) (v)=\phi_n'(a)\phi_n'(bc)(v)=\phi_n'(abc)(v)$,
(\ref{3nov15}) follows. \qed
\vskip 0.2in
Let $\rho:R\to\mat_{n\times n}(D)$ be an infinite
dimensional representation, where $D$ is a skew field
over the base field $K$. Clearly, there exists a countable dimensional
subskew field $E\subset D$ such that $\Ima(\rho)\subset \mat_{n\times n}(E)$.
Since $\mat_{n\times n}(E)$ is sofic if and only if $E$ is sofic and
$\mat_{n\times n}(E)$ has a unique rank, by Proposition \ref{P2}, we have that if $E$ is sofic, then $\rk_\rho$
is in the closure of $\Fin(R)$.
It is an open question, whether there are non-sofic skew fields or not.
In \cite{Eleklinear} it was proved that all the amenable skew fields
(we discuss them in Section \ref{skew}) and the free skew field (we will
prove that it is a non-amenable skew field in Section \ref{skew}) are sofic.
Hence for these skew fields and homomorphisms  $\rho:R\to\mat_{n\times n}(E)$,
$\rk_\rho$ is always in the closure of $\Fin(R)$.

\section{Amenable elements in 
the rank spectrum}
As we have seen in the previous section, one can construct convergent
sequences of finite dimensional representations of a finite dimensional
algebra $R$ by mapping $R$ into a sofic algebra $\cA$. Now we investigate that
using this construction, how can we obtain amenable elements (see Definition \ref{def1.1})in the
rank spectrum of $R$. 
\noindent
So, let $\cA$ be a countable dimensional algebra over the finite field $K$,
with basis $\{1=e_1,e_2,\dots\}$
and let $\{\phi_n:\cA\to \mat_{m_n\times m_n}(K)\}^\infty_{n=1}$ be unital
linear maps forming a convergent, sofic representation system. We
call the family $\{\phi_n\}^\infty_{n=1}$ hyperfinite (see \cite{Eleklinear})
if for any
$\e>0$, there exists $L_\e>0$ such that for any $n\geq 1$, we
have independent $K$-linear subspaces
$N^1_n,N^2_n,\dots, N_n^{k(n)}\subset K^{m_n}$ satisfying the following three 
conditions:
\begin{itemize}
\item For any $1\leq i \leq k(n)$, $\dim_K(N^i_n)\leq L_\e$
\item$\frac {\sum_{i=1}^{k(n)} \dim_K(N^i_n)}{m_n}\geq 1-\e$
\item The subspaces $\bigvee_{a\in \Span\{e_1,e_2,\dots,e_s\}}\phi_n(a)(N^i_n) $ are independent and for any $1\leq i\leq k(n)$,
$\frac{\dim_K\left(\bigvee_{a\in \Span\{e_1,e_2,\dots,e_s\}}\phi_n(a)(N^i_n)\right)}{\dim_K(N^i_n)}\leq 1+\e\,,$
where $s$ is the integer
part of $1/\e$.
\end{itemize}
\noindent
So, as opposed to the case of finite dimensional algebras, the pieces
 $N^i_n$ are not exactly $\cA$-modules, only in an approximate sense.
First, we prove a technical lemma that will make the proof of some of our results easier. It states that part of the third condition is
not necessary to check in order to prove that a certain sequence is hyperfinite. 
\begin{lemma} \label{technical}
Let $\{\phi_n:\cA\to \mat_{m_n\times m_n}(K)\}^\infty_{n=1}$ be a convergent sofic representation sequence.
Suppose that for any $\e>0$, there exists
$L_\e>0$ such that for any $n\geq 1$, we
have independent $K$-linear subspaces
$Z^1_n,Z^2_n,\dots, Z_n^{k(n)}\subset K^{m_n}$ satisfying the following three 
conditions:
\begin{itemize}
\item For any $1\leq i \leq k(n)$, $\dim_K(Z^i_n)\leq L_\e$.
\item$\frac {\sum_{i=1}^{k(n)} \dim_K(Z^i_n)}{m_n}\geq 1-\e$.
\item For any $1\leq i\leq k(n)$,
$\frac{\dim_K\left(\bigvee_{a\in \Span\{e_1,e_2,\dots,e_s\}}\phi_n(a)(Z^i_n)\right)}{\dim_K(Z^i_n)}\leq 1+\e\,,$
where $s$ is the integer
part of $1/\e$.
\end{itemize}
\noindent
Then $\{\phi_n\}$ is hyperfinite.
\end{lemma}
\proof We start with a simple linear algebra lemma.
\begin{lemma}\label{linear}
Let $V\subset Z$ be finite dimensional $K$-spaces and $T,S\subset \Endo_K(Z)$. 
Suppose that $\dim_K(T(V)+V)\leq (1+\e_T)\dim_K(V)$,  
$\dim_K(S(V)+V)\leq (1+\e_S)\dim_K(V)$,
then
\begin{itemize}
\item There exists a linear subspace $W\subset V$ such that $T(W)\subset V$ 
and $\dim_K(W)\geq  (1-e_T)\dim_K(V)\,.$
\item $\dim_K(TS(V)+V)\leq (1+e_T+e_S)\dim_K(V)$.
\end{itemize}
\end{lemma}
\proof Let $m_T:V\to \frac{T(V)+V}{V}$ be the natural quotient map. Since \\
$\dim_K(\Ima(m_T))\leq \e_T\dim_K(V)$, we have that
$$\dim_K(\Ker(m_T))\geq (1-\e_T)\dim_K(V)\,.$$
\noindent
On the other hand, $T(\Ker(m_T))\subset V$. For the second part, 
let $S(V)+V=Q\oplus V$ for
some complementing space $Q$, then $TS(V)+V\subset T(Q)+T(V)+V$, hence
$$\dim_K(TS(V)+V)\leq (1+\e_T+\e_S)\dim_K(V)\quad\qed\,.$$
\begin{lemma}
For any $n\geq 1$ and $1\leq i\leq k(n)$,
there exists a subspace $Q^i_n\subset Z^i_n$ so that
\begin{itemize}
\item
$$\dim_K (Q^i_n)\geq (1-\sqrt{\e})\dim_K (Z^i_n)\,.$$
\item
For any $a\in\Span\{e_1,e_2,\dots,e_{s(\e)}\}$
$$\phi_n(a)(Q^i_n)\subset Z^i_n\,,$$
where $s(\e)$ is given in such a way that $|\Span\{e_1,e_2,\dots,e_{s(\e)}\}|\leq \frac{1}{\sqrt{\e}}$.
\end{itemize}
\end{lemma}
\noindent
\proof
For $a\in\Span\{e_1,e_2,\dots,e_{s(\e)}\}$, let $W_a\subset Z^i_n$ (Lemma \ref{linear})  such that
\begin{itemize}
\item $\dim_K(W_a)\geq (1-\e)\dim_K(Z^i_n)\,.$
\item $\phi_n(a)(W_a)\subset Z^i_n\,.$
\end{itemize}
\noindent
Let $$Q^i_n=\cap_{a\in\Span\{e_1,e_2,\dots,e_{s(\e)}\}} W_a\,.$$
\noindent
Then,
$$\dim_K(Q^i_n)\geq (1-\sqrt{\e})\dim_K(Z^i_n)$$
and
for any $a\in\Span\{e_1,e_2,\dots,e_{s(\e)}\}$
$$\phi_n(a)(Q^i_n)\subset Z^i_n\quad\qed$$
\noindent
The existence of the systems $\{Q^i_n\}^{k(n)}_{i=1}$ clearly implies the statement in Lemma \ref{technical} \qed
\vskip 0.1in
\noindent
If $\{\phi_n:\cA\to \mat_{m_n\times m_n}(K)\}^\infty_{n=1}$ is a hyperfinite
sofic representation system, then the associated sofic representation
system
$$\{\hat{\phi}_n:\mat_{l\times l}(\cA)\to \mat_{m_n l\times m_n l}(K)\}$$
is still hyperfinite. Indeed, the spaces $M^i_n:=(N^i_n)^l\subset K^{m_nl}$ will
be the independent subspaces satisfying the approximate module condition.
Now, let $\rho: R\to\mat_{l\times l}(\cA)$ be a unital homomorphism from
a finite dimensional algebra $R$. By Theorem \ref{P2}, we have
a sequence of subspaces $V_n\subset K^{m_nl}$ so that
\begin{itemize}
\item $\hat{\phi}_n\circ \rho(V_n)\subset V_n$ (that is $V_n$ is an $R$-module).
\item $\lim_{n\to\infty} \frac{\dim_K(V_n)}{m_n l}=1\,,$
\end{itemize}
\noindent
representing the element $\rk_\rho$ in the rank spectrum.
\begin{propo} \label{transfer}
$\rk_\rho$ is an amenable element of the rank spectrum, that is, 
$\{V_n\}^\infty_{n=1}$ is a hyperfinite sequence of $R$-modules.
\end{propo}
\proof The first step is to show that the
sequence $\{V_n\}^\infty_{n=1}$ is hyperfinite
in the weaker sense, described in Lemma \ref{technical}. As opposed to the graph theoretical case, it is
not a priori obvious. The reason is that although
the subspaces $\{V_n\}^\infty_{n=1}$ are
almost as big as the total spaces $K^{m_n l}$ it is not true that they
contain all (or even one single copy) of the pieces $M^i_n$. The following
lemma resolves this problem. First, fix a basis $\{f_i\}^\infty_{i=1}$ for the algebra $\mat_{l\times l}(A)$.
By taking a subsequence, we can suppose that
the subspaces $\bigvee_{a\in \Span\{f_1,f_2,\dots,f_s\}}\hat{\phi}_n(a)(M^i_n) $ are independent, where $s$ is the integer part of $1/\e$ and
for any $1\leq i\leq k(n)$,
$$\frac{\dim_K\left(\bigvee_{a\in \Span\{f_1,f_2,\dots,f_s\}}\hat{\phi}_n(a)(M^i_n)\right)}{\dim_K(M^i_n)}\leq 1+\e\,,$$
\begin{lemma}
For any $\e>0$ and large enough $n\geq 1$, we have independent
subspaces $Q^1_n,Q^2_n,\dots Q^{l(n)}_n$ in $V_n$ such that
\begin{itemize}
\item For any $1\leq i \leq l(n)$, $\dim_K(Q^i_n)\leq L_\e\,.$
\item$\frac {\sum_{i=1}^{l(n)} \dim_K(Q^i_n)}{m_n}\geq 1-2\e\,.$
\item For any $1\leq i\leq l(n)$,
$\frac{\dim_K\left(\bigvee_{a\in \Span\{f_1,f_2,\dots,f_s\}}\hat{\phi}_n(a)(Q^i_n)\right)}
{\dim_K(Q^i_n)}\leq 1+\e\,,$
where $s$ is the integer part of $1/\e$.
\end{itemize}
\end{lemma}
\proof
We call the subspaces $M^i_n$ and $M^j_n$ equivalent if there exists a linear bijection
$\alpha_{ij}:K^{m_nl}\to K^{m_nl}$ so that
\begin{itemize}
\item
$\alpha_{ij}$ bijectively maps $M^i_n$ into $M^j_n$.
\item
For any $a\in  \Span\{f_1,f_2,\dots,f_s\}$ and $v\in M^i_n$
$$\alpha_{ij}(\hat{\phi}_n(a)(v))=\hat{\phi}_n(a)(\alpha_{ij}(v))\,.$$
\end{itemize}
Notice that by our finiteness conditions, there exists a constant $C>0$ such
that for any $n\geq 1$, the number of equivalence classes is at most $C$.
Now, let $\mu=\{M^{i_1}_n, M^{i_2}_n,\dots,M^{i_{\mu(t)}}_n\}$ be an 
equivalence class
and let $\alpha_j:K^{m_nl}\to K^{m_nl}$ be the 
corresponding bijections mapping $M^{i_1}_n$ into $M^{i_j}_n$.
For $\hat{\lambda}=\{\lambda_1,\lambda_2,\dots,\lambda_{i_{\mu(t)}}\}
\in K^{\mu(t)}$ define
$$M^{\hat{\lambda}}_n:=(\sum_{j=1}^{\mu(t)} \lambda_j\alpha_j)(M^{i_1}_n)\,.$$
\noindent
Then
$$\dim_K\left(\bigvee_{a\in \Span\{f_1,f_2,\dots,f_s\}}\hat{\phi}_n(a)
(M^{\hat{\lambda}}_n)\right)
\leq (1+\e)\dim_K(M^{\hat{\lambda}}_n)\,.$$
\noindent
Also, if $\hat{\lambda}_1, \hat{\lambda}_2,\dots,
\hat{\lambda}_q$ are independent vectors in $K^{\mu(t)}$, then
the corresponding subspaces $M^{\hat{\lambda}_1}_n,
M^{\hat{\lambda}_2}_n,\dots,M^{\hat{\lambda}_q}_n$ are independent as well.
We call an equivalence class $\mu$ {\it large} if
$$\mu(t)\dim_K
(M^{\mu}_n)\geq \frac{(1-\e)\e}{100 C} m_nl $$
\noindent
holds, where $M^{\mu}_n$ is representing the class $\mu$. It is not hard
to check that
\begin{equation}\label{large}
\dim_K(\bigoplus_{\mbox{$\mu$ is large}} \oplus_{i=1}^{\mu(t)} M^{i,\mu}_n)
\geq (1-\frac{3}{2}\e) m_n l\,,
\end{equation}
holds, where $\{M^{1,\mu}_n,M^{2,\mu}_n,\dots,M^{\mu(t),\mu}_n\}$ are the
subspaces in the class $\mu$.
Now let $\mu_n$ be a large class for some $n\geq 1$ and $M_n^{i_1}$ be 
its representative.
For $v\in M^{i_1}_n$, let $H^{\mu_n}_v\subset K^{\mu_n(t)}$ be the set of vectors
$\hat{\lambda}$ such that $\sum^{\mu_n(t)}_{j=1} \lambda_j \alpha_j(v)\in V_n.$
Since $\lim_{n\to\infty}\frac{\dim_K(V_n)}{m_nl}=1,$ we get
\begin{equation}\label{masiklarge}
\lim_{n\to\infty}\frac{\dim_K(\cap_v H^{\mu_n}_v)}{\mu_n(t)}=1\,.
\end{equation}
\noindent
For each large class $\mu_n$, when $n$ is large enough we choose
a basis $\hat{\lambda}_1,\hat{\lambda}_2,\dots,\hat{\lambda}_p$ in
$\cap_v H^{\mu_n}_v$ and set $Q^n_{i,\mu_n}:=M^{\hat{\lambda}_i}_n\,.$ Then, by
(\ref{large}) and (\ref{masiklarge}) the subspaces $Q^n_{i,\mu_n}$
will satisfy the conditions of our proposition, provided that $n$ is
large enough. \qed  
\vskip 0.2in
\noindent
So, we have a convergent sequence of representations
$$\{\psi_n:R\to\Endo_K(V_n)\}^\infty_{n=1}$$
that are hyperfinite in the weaker sense, that is for any $\e>0$
we have $L_\e>0$ and for large enough $n$ independent subspaces
$Q^1_n,Q^2_n,\dots,Q^{l(n)}\subset V_n$ such that
\begin{itemize}
\item For any $1\leq i \leq l(n)$, $\dim_K(Q^i_n)\leq L_\e$.
\item $\sum_{i=1}^{l(n)} \dim_K(Q^i_n) \geq (1-2\e)\dim_K(V_n)\,.$
\item  For any $1\leq i \leq l(n)$,
$$\dim_K\left(\bigvee_{r\in R} \psi_n(r)(Q^i_n)\right)\leq
(1+\e)\dim_K(Q^n_i)\,.$$
\end{itemize}
\noindent
Note that we used the fact that if $s$ is large then
$\Span\{f_1,f_2,\dots,f_s\}$ contains $\rho(R)$.
Our proposition immediately follows from the following lemma.
\begin{lemma}
Let $\tau:R\to \Endo_K(V)$ be a unital representation and $L\subset V$
be a subspace such that
$$\dim_K\left(\bigvee_{r\in R} \tau(r)(L)\right)\leq (1+\e)\dim_K(L)\,.$$
\noindent
Then there exists a subspace $L'\subset L$, $\frac{\dim_K(L')}{\dim_K(L)}\geq 1-|R|\e$
such that $L'$ is an $R$-module (that is for any $r\in R$, $\tau(r)L'\subset L'$).
\end{lemma}
\noindent
\proof
By Lemma \ref{linear}, for any $r\in R$, there
exists a subspace  $L_r\subset L$ such that
$$\dim_K(L_r)\geq (1-\e)\dim_K(L)$$
and
$\tau(r)L_r\subset L.$
Then $\dim_K(\cap_{r\in R} L_r)\geq (1-|R|\e)\dim_K(L)$
and for
any $r\in R$, $\tau(r)(\cap_{r\in R} L_r)\subset L$.
Hence the subspace $L'=\bigvee_{r\in R}( \tau(r)(\cap_{r\in R} L_r))$ satisfies
the condition of our lemma. \qed

\begin{propo}\label{poly}
Let $K[X]$ be the polynomial algebra over $K$. Then the family 
of all the finite dimensional representations $\kmod$ are hyperfinite.
\end{propo}
\begin{corol}
If $\phi:R\to \mat_{l\times l}(K[X])$ is a homomorphism, then
$\phi^*(\kmod)$ is a hyperfinite family in $\rmod$.
\end{corol}
\proof (of Proposition \ref{poly})
It is enough to show that the class
of indecomposable finite dimensional $K[X]$-modules is hyperfinite. By Jordan's theorem
indecomposable elements in $\kmod$ are in the form of $K[X]/f^nK[X]$, where $n\geq 1$
and $f\in K[X]$ is an irreducible, monic polynomial. The module structure of
$M_{f^n}=K[X]/f^nK[X]$ is given the following way,
\begin{itemize}
\item $M_{f^n}=\Span\{1,x,x^2,\dots,x^{\deg(f)n-1}\}$
\item $X\cdot t^i=t^{i+1}$, if $1\leq i \leq \deg(f)n-2$
\item $X\cdot t^{\deg(f)n}-1=t^{\deg(f)n}-f^n(t)\,.$
\end{itemize} 
\noindent
Let $\e>0$ and $k\geq 1/\e$. We will show that for any module $M_{f^n}$, we have
independent subspaces $N^1_{f^n}, N^2_{f^n},\dots, N^{k(f^n)}_{f^n}\subset M_{f^n}$ such that
\begin{itemize}
\item For any $1\leq i \leq k(f^n)$, $\dim_K(N^i_{f^n})\leq 2k^2$
\item  For any $1\leq i \leq k(f^n)$, 
$$\frac{\dim_K(X(N^i_{f^n})+N^i_{f^n})}{\dim_K(N^i_{f^n})}\leq 1+\e$$
\item $\sum_{i=1}^{k(f^n)} \dim_K(N^i_{f^n})\geq (1-\e)\dim_K(M_{f^n})\,.$
\end{itemize} 
\noindent
If $\deg(f)n\leq 2k^2$, let $k(f^n)=1$ and $N^1_{f^n}=M_{f^n}.$ 
If $\deg(f)n > 2k^2$, let C be the integer such that 
$$Ck-1 <\deg(f)n\leq  (C+1)k-1\,$$
and for $1\leq i \leq C+1$ let
$$N^i_{f^n}:=\Span\{t^{(i-1)k}, t^{(i-1)k}+1,\dots, t^{ik-1}\}\,.$$
\noindent
Then $\dim_K(X(N^i_{f^n})+N^i_{f^n})= \dim_K(N^i_{f^n})+1$ and all the three
inequalities above are satisfied. \qed
\section{Amenable and non-amenable skew fields}
\label{skew}
In this section we recall some basic definitions and results about amenable
and non-amenable skew fields. We will confirm Conjecture \ref{con2} for
amenable algebras $\mat_{l\times l}(D)$, when $D$ is an amenable skew field.
We also prove the non-amenability of the free skew 
field $D^r_K$. This result is hopefully interesting on its own right, 
nevertheless it will be important for us in the proof
of Theorem \ref{kronecker}.

\noindent
 Let $\cA$ be a countable dimensional algebra over
the base field $K$. We say that the algebra is (left) amenable (\cite{Gromov},
\cite{Elekskew}) if there exists a sequence of finite dimensional 
$K$-linear subspaces $\{W_n\}^\infty_{n=1}\subset \cA$ such that
for each $a\in \cA$
$$\lim_{n\to\infty} \frac{\dim_K(W_n+aW_n)}{\dim_K(W_n)}=1\,.$$
\noindent
Let us recall some basic facts on amenable algebras.
\begin{itemize}
\item If $\cA$ is an amenable domain, then it has the Ore property and its 
classical skew field of quotients is amenable as well (Proposition 2.1, 
\cite{Elekskew}).
\item If $\cA=K\Gamma$ is a group algebra, then $\cA$ is amenable if
and only if $\Gamma$ is an amenable group \cite{Elekskew}, \cite{Bartholdi}.
\item If $E\subset D$ are skew fields over the base field $K$, and
$E$ is non-amenable, then so is $D$ (Theorem 1. \cite{Elekskew})
\item If $D$ and $E$ are amenable skew fields and $D\otimes_K E$ is
a domain, then $D\otimes_K E$ is amenable (Proposition 2.4 \cite{Elekskew}).
\end{itemize}
Now let $KF_r$ be the free algebra $K\langle x_1,x_2,\dots,x_t\rangle$ on 
$r$ indeterminates. Let $\fmod$ be the set of all finite dimensional modules over
$KF_r$. So, any element $M\in\fmod$ defines a Sylvester rank function on $KF_r$
and the convergence of such representations is well-defined (see \cite{Eleklinear}) exactly
the same way as for finite dimensional algebras. Note that a convergent sequence
of representations can always be viewed as a sofic approximation of the
algebra $KF_r/I$, where $I$ is the ideal of elements $a\in KF_r$ for which 
$\lim_{n\to\infty} \rk_{M_n}(a)=0.$

\noindent
Let $\cM=\{M_n\}^\infty_{n=1}\subset \fmod$ be a convergent sequence
of modules and $\phi_\cM:KF_r\to\mat^\cM_\omega$ be the ultraproduct
representation as in the previous sections (we suppose
that $K$ is a finite field). Then the division closure of
the image of $KF_r$ in $\mat\cM_\omega$ is a skew field if and only if
the rank function $\rk_\omega$ is integer-valued (Theorem 2. \cite{Eleklinear}).

\noindent
Furthermore, $D$ is an amenable skew field if and only if the convergent
sequence $\cM$ is hyperfinite (Theorem 3. \cite{Eleklinear}). According to 
Proposition 7.1 \cite{Eleklinear} the free skew field $D^r_K$ on $r$
generators over the base field $K$ is sofic
(for the free skew field see e.g.\cite{Cohn}.) That is, there exists
a convergent
sequence $\{M_n\}^\infty_{n=1}\subset \fmod$  , for which the division closure in $\mat^\cM_\omega$ is
the free skew field.
In \cite{Elekskew} using operator algebraic methods we proved that the
free field over the complex numbers is non-amenable. Now we prove the 
following theorem.
\begin{theorem}\label{free}
The free skew field on $r>1$ generators over an arbitrary base field $K$ is
non-amenable.
\end{theorem}
\proof
First, observe that for any base field $K$, the free algebra
$KF_r$ is isomorphic to its opposition algebra
$KF_r^{op}$, hence the free skew field $D=D^r_K$ is
isomorphic to its opposite algebra $D^{op}$. Recall that
$D\otimes_K D$ is a domain (Theorem 3.1, \cite{Cohncoll}). Hence,
it is enough to prove that there exists an embedding of
$D\otimes D^{op}$-modules
\begin{equation} \label{egy1}
\pi:\left( D\otimes D^{op}\right)^2\to D\otimes D^{op}\,.
\end{equation}
\noindent
Indeed, if $D$ is amenable, then by the properties listed above
$D\otimes D^{op}$ is an amenable domain and therefore 
$D\otimes D^{op}$ has the Ore property. Thus two non-trivial
left ideals of $D\otimes D^{op}$ have non-zero intersection, in
contradiction with the existence of the embedding (\ref{egy1}).
In order to construct the embedding, we follow the lines of
Lemma 4.21 in \cite{BCL} (I am indebted to Andrey Lazarev for
calling my attention to this paper). First, consider the two-terms
free resolution of $\KR$ as a bimodule over itself 
$$0\to\KR\otimes_K[x_1,x_2]\otimes_K \KR\to\KR\otimes_K\KR\to
\KR\to 0$$
\noindent
as in \cite{BCL}. Here $[x_1,x_2]$ denotes the 
2-dimensional vector space spanned by $x_1,x_2$.
 If we tensor the resolution above by $D$ on the right, we obtain an
exact sequence of $\KR-D$\,-bimodules
$$0\to\KR\otimes_K[x_1,x_2]\otimes_K D\to\KR\otimes_K D \to
D\to 0\,,
$$
\noindent
since $\Tor^{\KR}_1(\KR,D)=0$ by the definition of the functor $\Tor$.
Now let us tensor this exact sequence by $D$ on the left.
Notice that  $\Tor^{\KR}_1(D,D)=0$ (\cite{Scho} Theorems 4.7,4.8).
Hence, we obtain an embedding of $D-D$\,-bimodules
$$i:D\otimes_K[x_1,x_2]\otimes_K D\to D\otimes D\,.$$
That is we have the embedding of left $D\otimes D^{op}$-modules
$$j:(D\otimes D^{op})^2\to D\otimes D^{op}$$
we sought for. This finishes the proof of our theorem. \qed
\begin{theorem}\label{fontostetel}
Let $\rho:R\to\mat_{l\times l}(D)$ be an infinite dimensional
representation of a finite dimensional algebra $R$ over the finite
field $K$, where $D$ is a countable dimensional amenable skew field. Then the
associated element $\rk_\rho\in \Rank(R)$ is amenable.
\end{theorem}
\proof
Let $\{d_1,d_2,\dots,d_r\}$ be the set of all entries in the matrices $\{\rho(s)\}_{s\in R}.$
Let 
$\pi:KF_r\to D$ be the unital homomorphism mapping $x_i$ to $d_i$.
Let $E\subset D$ be the sub skew field generated by $\{ d_1, d_2,\dots, d_r\}$ (that is
the division closure of $\Ima(\pi))$.
From now on, we can suppose that $\rho$ maps $R$ into $\mat_{l\times l}(E)$ and $\pi$
maps $KF_r$ into $E$.
Let $S:=KF_r/\Ker(\pi)$, and $\theta:KF_r\to S$ and $\hat{\pi}:S\to E$ be the natural
quotient maps, so $\pi=\hat{\pi}\circ \theta$.
Then, the homomorphism $\zeta:R\to\mat_{l\times l}(S)$ can be defined in a unique way to
satisfy $\hat{\pi}\circ \zeta=\rho$ (note that $\hat{\pi}$ is injective).
Finally, pick a lifting $\mu:S\to KF_r$. that is $\mu$ is an injective, unital, linear
map such that $\theta\circ \mu=\Id_S$. Then, $\pi\circ \mu=\hat{\pi}$.

\noindent
Since $E$ is amenable, it is sofic. Hence, we have a sofic representation sequence
$\{\phi_n:E\to \mat_{m_n\times m_n}(K)\}^\infty_{n=1}$. Therefore, we have
the sofic representation sequence
$$\{\phi_n\circ \hat{\pi}:S \to\mat_{m_n\times m_n}(K)\}^\infty_{n=1}$$
\noindent
that factors through the sofic representation system
$$\{\phi_n\circ \pi:KF_r \to\mat_{m_n\times m_n}(K)\}^\infty_{n=1}\,.$$
\noindent
By Proposition 11.1 \cite{Eleklinear}, $\{\phi_n\circ \pi\}^\infty_{n=1}$ is hyperfinite
hence $\{\phi_n\circ \hat{\pi}\}^\infty_{n=1}$ is a hyperfinite sofic representation.
Therefore by Proposition \ref{transfer}, $\rk_\zeta=\rk_\rho$ is an amenable element
of the rank spectrum. \qed

\section{The wild Kronecker algebras are of non-amenable representation types}
First recall the notion of Kronecker quiver algebras
and the wildness
phenomenon  (see e.g. \cite{Barot}). Let $r\geq 3$ and $Q_r$ be the  finite dimensional algebra
with basis $p_1,p_2,\{e_i\}^r_{i=1}$, where
\begin{itemize}
\item $p_1^2=p_1, p_2^2=p_2, p_1+p_2=1, p_1p_2=p_2p_1=0$.
\item For any $1\leq i,j\leq r$, $e_ip_1=e_i$, $p_1e_i=0$, 
$e_ip_2=0$, $p_2e_i=e_i$, $e_ie_j=0$.
\end{itemize}
\noindent
The algebra $Q_r$ defined the way above is the Kronecker quiver algebra
of index $r$.
The wildness of $Q_r$ means that the classification of finitely generated
indecomposable modules over $Q_r$ is as complicated as the classification of 
finitely
generated indecomposable modules over noncommutative free algebras. 
To make it precise, there exists a representation
$\pi_r:Q_r\to \mat_{2\times 2}(KF_{r-1})$ such that
\begin{itemize}
\item The associated functor $\pi_r^*:\flmod\to \qmod$ is injective.
\item If $M\in\flmod$ is indecomposable, then $\pi_r^*(M)$ is 
indecomposable as well.
\end{itemize}
\noindent
The following representation satisfies the conditions above:
$\pi_r(p_1)=\begin{bmatrix} 1 & 0 \\ 0 & 0 \end{bmatrix},$
$\pi_r(p_2)=\begin{bmatrix} 0 & 0 \\ 0 & 1 \end{bmatrix},$
$\pi_r(e_1)=\begin{bmatrix} 0 & 0 \\ 1  & 0 \end{bmatrix},$
$\pi_r(e_j)=\begin{bmatrix} 0 & 0 \\ x_{j-1}  & 0 \end{bmatrix}$  if
$j>1$.
\begin{propo}\label{ingredient}
Let $\pi_r:Q_r\to \mat_{2\times 2}(KF_{r-1})$ the representation as above. Then
the family of modules $\{\pi_r^*(M_n)\}^\infty_{n=1}\subset \qmod$ is
hyperfinite if and only if  $\{M_n\}^\infty_{n=1}\subset \flmod$ is
hyperfinite.
\end{propo}
\proof The ``if'' part follows from Proposition \ref{transfer}. For the converse,
let us suppose that  $\{\pi_r^*(M_n)\}^\infty_{n=1}$ is a hyperfinite family.
That is, for any $\e>0$, we have $L_\e>0$ and for any $n\geq 1$ we have
independent $Q_r$-submodules $N^1_n,N^2_n,\dots, N^{k(n)}_n\subset \pi_r^*(M_n)$
such that
\begin{itemize}
\item $\dim_K(N^i_n)\leq L_\e$
\item $\sum_{i=1}^{k(n)} \dim_K(N^i_n)\geq (1-\e)2\dim_K(M_n)$.
Note that as a vector space $\pi_r^*(M_n)$ is isomorphic to $M_n\oplus M_n$ and $p_i$ acts
on $\pi_r^*(M_n)$ as the projection onto the $i$-th coordinate.
\end{itemize}
\noindent
Let us consider the subspaces $p_1N^i_n\subset M_n\oplus 0$, 
$p_2N^i_n\subset 0\oplus M_n$. We use the notation $[p_1 N^i_n]$ resp.
 $[p_2 N^i_n]$
for the subspaces in $M_n$ that are the projections of the spaces
above onto the first
resp. second component of $M_n\oplus M_n$. Then, for any $n\geq 1$,
$1\leq i \leq k(n)$
\begin{itemize}
\item $\pbe\subset \pbk$
\item $x_j\pbe\subset \pbk$ if $1\leq j\leq r-1$.
\end{itemize}
\noindent
Also, the subspaces $\{\pbe\}^{k(n)}_{i=1}$ are independent.
Since
$$\sum^{k(n)}_{i=1} (\dim_K\pbe +\dim_K\pbk)\geq (1-\e) 2\dim_K(M_n)\,,$$
\noindent
we have the inequality
$$\sum^{k(n)}_{i=1} \dim_K\pbe\geq (1-2\e) \dim_K(M_n)\,.$$
\begin{lemma}
Let  $C>0$, $\{a_i\}^{k(n)}_{i=1}$,  $\{b_i\}^{k(n)}_{i=1}$, be positive
real numbers satisfying the following inequalities.
\begin{itemize}
\item $\sum^{k(n)}_{i=1} a_i\geq (1-2\e)C$
\item $\sum^{k(n)}_{i=1} b_i\leq C$
\item $a_i\leq b_i$.
\end{itemize}
\noindent
Let $S=\{i\,\mid \frac{b_i}{a_i}< 1+\sqrt{\e}\}.$
Then $\sum_{i\neq S} a_i\leq 2\sqrt{\e} C\,.$
\end{lemma}
\proof
We have
$$ C\geq \sum^{k(n)}_{i=1} b_i\geq \sum_{i\notin S} (1+\sqrt{\e})a_i +
\sum_{i\in S}a_i\geq $$
$$\geq (1-2\e)C+ (\sum_{i\notin S}a_i)\sqrt{\e}\,.$$
Therefore $\sum_{i\neq S} a_i\leq 2\sqrt{\e} C\,.$ \qed
\vskip 0.2in
\noindent
Let us apply the lemma above, for $C=\dim_K(M_n)$, $a_i=[p_1 N^n_i]$,
$b_i=[p_2N^n_i]\,.$ We get that
\begin{itemize}
\item $\sum_{i\notin S} \dim_K([p_1N_i^n])\leq 2\sqrt{\e}\dim_K(M_n)\,.$
\item If $i \in S$, then $ \dim_K\left([p_1N_i^n]\oplus \bigoplus_{j=1}^{r-1}x_j[p_1N_i^n]\right)
\leq (1+\sqrt{\e})
\dim_K([p_1N_i^n])$ .
\item $\dim_K([p_1N_i^n])\leq L_\e$.
\end{itemize}
Hence by definition, $\{M_n\}^\infty_{n=1}$ is a hyperfinite family of
$KF_{r-1}$-modules. \qed
\vskip 0.2in
\begin{theorem} \label{kronecker}
The wild Kronecker quiver algebras are of non-amenable representation type.
\end{theorem}
\proof
Let $r\geq 2$ and  $\cM=\{M_n\}^\infty_{n=1}\subset \fmod$ be
a convergent sequence of $KF_r$-modules such that for all $KF_r$-matrices $A$
$$\lim_{n\to\infty} \rk_{M_n}(A)\in \Z\,.$$
\noindent
Let $\phi_\omega:KF_r\to\mat^\cM_\omega$ be the limit representation of
the sequence $\{M_n\}^\infty_{n=1}$ Then by Proposition 8.1 \cite{Eleklinear},
the division closure of $\phi_\omega(KF_r)$ in the ring $\mat^\cM_\omega$ is
a skew field. We also know that
\begin{itemize}
\item
The division closure above is an amenable skew field if and only if
$\{M_n\}^\infty_{n=1}$ is hyperfinite 
(Proposition 11.1 and 12.1 \cite{Eleklinear}).
\item
There exists a convergent sequence of modules  $\{M_n\}^\infty_{n=1}$
such that the division closure is the free skew field on $r$-generators
with base field $K$ (Proposition 7.1 \cite{Eleklinear}).
\end{itemize}

\noindent
Hence by Proposition \ref{ingredient} and
Theorem \ref{free}, there exist non-hyperfinite families of modules for
the algebra $Q_r$, $r\geq 3$. \qed

\section{Hyperfinite families and the continuous algebra of 
von Neumann} \label{hypersection}

The goal of this section is to prove Theorem \ref{embedding}. In fact, we
will prove a stronger interpolation theorem.  Let $K$ be 
a finite field, $R$ be a finite dimensional algebra over $K$ and
let $\cM=\{M_n\}^\infty_{n=1}\subset \rmod$ be a hyperfinite convergent sequence
of modules and $\omega$ be a non-principal
ultrafilter. 
\begin{theorem} \label{tetel2}
Let $\rho_\omega:R\to \mat^\cM_\omega$ be the ultraproduct representation
 corresponding
to the sequence of modules above. Then there exists a subalgebra
$\Ima(\rho_\omega)\subset S \subset \mat^\cM_\omega$ such that  $S\cong M_K$.
\end{theorem}
\noindent
According to the theorem, when we pick a countable dimensional subalgebra
from the ultraproduct, we can always land inside $M_K$.
\proof The proof will be given in a series of lemmas and propositions.
\begin{lemma} \label{aug22}
Let $a_1,a_2,\dots, a_r$ be positive integers, $\e>0$ and $U\in\omega$. 
Let us have an integer $m_i$ for each $i\in U$ (we suppose that $m_i$
tends to infinity as $i$ tends to infinity).
Suppose that for any $i\in U$ and $1\leq j\leq r$
we have a non-negative integer $c_{ij}$ such that
$$(1-\e) m_i\leq \sum^r_{j=1} c_{ij} a_j\leq m_i\,.$$
\noindent
Then we have a subset $V\subset U$, $V\in\omega$, and non-negative integers
$\{p_j\}^r_{j=1}$, $\{l_k\}^\infty_{k=1}$ such that for any
$i\in V$
$$(1-2\e) m_{i}\leq l_i (\sum^r_{j=1} p_j a_j)\leq m_{i}\,.$$
\end{lemma}
\proof
For large enough $i\in U$, let
$$\frac{\e m_i}{2r(\max_{1\leq j\leq r} a_j)}<l_i< 
\frac{\e m_i}{r(\max_{1\leq j\leq r} a_j)}$$
be an integer. Then, for any $1\leq j\leq r$, we have
$c_{ij}=l_id_{ij}+g_{ij}$, where $0\leq g_{ij}<l_i$ and 
$$0\leq d_{ij}\leq \frac{2r\max_{1\leq j\leq r} a_j}{\e}\,.$$

\noindent
Since
$$\sum^r_{j=1}l_i a_j\leq \frac{\e m_i} {r(\max_{1\leq j\leq r} a_j)}
\sum^r_{j=1} a_j\leq \e m_i\,,$$
we have that $\sum^r_{j=1} g_{ij}a_j\leq \e m_i$ that is
$$(1-2\e)m_i\leq l_i\sum^r_{j=1}d_{ij}a_j\leq m_i\,.$$
\noindent
Since the set $\{d_{ij}\}_{i\in U, 1\leq j\leq r}$ is bounded, there are 
non-negative integers $p_1,p_2,\dots,p_r$ such that
$$V=\{i\mid\,d_{ij}=p_j,\,\mbox{for any $1\leq j \leq r$}\}\in\omega\,.$$
Therefore if $i\in V$, we have
$$(1-2\e) m_{i}\leq l_i (\sum^r_{j=1} p_j a_j)\leq m_{i}\,.\quad\qed$$
\begin{defin}
The module $A\in\rmod$ $\e$-tiles the module $B\in\rmod$, if there exists
$k\geq 1$ and $A^k\cong C\subset B$ such that $\dim_K(C)\geq (1-\e)\dim_K(B)$.
\end{defin}
\begin{lemma} \label{l52}
Let $U\in\omega$ be a subset of the naturals. Suppose that $\{M_i\}_{i\in U}$
is a hyperfinite family. Then for any $\epsilon$, there exists
$A\in\rmod$ and $V\subset U$, $V\in\omega$ such that
for each $i\in V$ the module $A$ $\e$-tiles $M_i$.
\end{lemma}
\proof
By hyperfiniteness, we have modules $A_1,A_2,\dots,A_r\in\rmod$ such that for
each $n\geq 1$ there exist constants $\{c_{nj}\}_{1\leq j \leq r}$ so that
$$\oplus^r_{j=1}A_j^{c_{nj}}\cong B_n\subset M_n\,$$
satisfying the inequality
$$(1-\frac{\e}{2})\dim_K(M_n)\leq \dim_K(B_n)\,.$$
By Lemma \ref{aug22}, there exist constants $\{p_j\}^r_{j=1}$,
a module $A=\oplus^r_{j=1} A_j^{p_j}$ and $V\subset U$, $V\in\omega$
so that $A$ $\e$-tiles $M_i$ if $i\in V$. \qed
\begin{lemma} \label{l53}
Suppose that the module $A\in\rmod$ $\e$-tiles all elements of the
sequence $\{B_n\}^\infty_{n=1}\subset\rmod$, where 
$\lim_{n\to\infty}\dim_K(B_n)=\infty$. Then for any $l\geq 1$, $A^l$
$2\e$-tiles $B_n$, if $n$ is large enough.
\end{lemma}
\noindent
(the proof is straightforward)
\begin{lemma} \label{l4aug22}
Let $0<\e<1$ and $A^k\cong C\subset B$ finitely generated
$R$-modules such that
$$(1-\e)\dim_K(B)\leq \dim_K(C)\,.$$
\noindent
Let $\delta\leq \frac{\e}{\dim_K^2(A)}$ and $B'\subset B$ is
a submodule such that
$$(1-\delta)\dim_K(B)\leq \dim_K(B')\,.$$
\noindent
Then, there exists $k'\geq 0$ and $A^{k'}\cong C'\subset B'$ such that
$$(1-2\e)\dim_K(B')\leq \dim_K(C')\,.$$
\end{lemma}
\proof
Let $b=\dim_K(B)$, that is, 
$$(1-\e)b\leq k\dim_K(A)\leq b\,.$$
\noindent
Let us write $A^k$ into the form of $A\otimes_K K^k$ and
for each $a\in A$ define the vector space
$$V_a:=\{v\,\mid a\otimes v\in B'\}\,.$$
\noindent
Since
$$\dim_K(a\otimes K^k)+\dim_K(B')=\dim_K\left((a\otimes K^k)\cap B'\right)+
\dim_K\left((a\otimes K^k)\vee B'\right)$$
we have that
$$k+(1-\delta)b\leq k+\dim_K(B')\leq \dim_K(V_a)+b\,.$$
\noindent
That is, $k-\delta b\leq\dim_K(V_a)\,.$ Hence,
$k-\dim_K(A)\delta b\leq \dim_K(\cap_{a\in A} V_a)\,.$
Observe that
$$C'=A\otimes_K(\cap_{a\in A} V_a)\subset B'$$
\noindent
and
$$(1-2\e)b\leq k \dim_K(A)-\delta b\dim_k^2(A)\leq \dim_K(C')\,.$$
\noindent
This finishes the proof of our lemma. \qed

\vskip 0.1in
\noindent
Now we put together the previous lemmas to have a single technical
proposition.
\begin{propo}\label{protech}
Let $\{M_n\}^\infty_{n=1}\subset \rmod$ be a convergent, hyperfinite
sequence. Then there exist modules
$\{N_i\}^\infty_{i=1}\subset\rmod$, subsets $\N\supseteq V_1\supseteq V_2\supseteq\dots,
V_i\in\omega$ and a monotonically 
decreasing sequence of real numbers $\{\alpha_n\}^\infty_{n=1}$, $\alpha_n\to 0$ such that
\begin{itemize}
\item For any $i$, $N_i$ $\alpha_i$-tiles $M_j$, whenever $j\in V_i$.
\item For any $i$, $N_i$ $\alpha_i$-tiles $N_{i+1}$.
\item For any $i$, $\dim_K(N_{i+1})$ is divisible by $\dim_K(N_i).$
\item For any $k,l\geq 1$ and matrix $A\in\mat_{k\times l}(R)$.
$$\lim_{n\to\infty} \rk_{M_n} (A)=\lim_{i\to\infty} \rk_{N_i} (A)\,.$$
\end{itemize}\end{propo}
\proof
We proceed by induction. Our inductional hypothesis goes as follows.
Suppose that $N_1,N_2,\dots N_l$ have already been constructed together
with the sets $V_1\supseteq V_2\supseteq\dots\supseteq V_l$ and $\alpha_1> \alpha_2 >\dots >\alpha_l$, $\alpha_i<\frac{1}{i},$
satisfying the following
conditions.
\begin{itemize}
\item For any $1\leq i\leq l-1$, $N_i$ $\al_i$-tiles $M_j$ if $j\in V_i$.
\item For any $1\leq i\leq l-1$, $N_i$ $\al_i$-tiles $N_{i+1}$.
\item $N_l$ $\al_l/2$-tiles $M_j$ if $j\in V_l$.
\item  For any $1\leq i\leq l-1$, $\dim_K(N_{i+1})$ is divisible by $\dim_K(N_i).$
\end{itemize}

\noindent
By Lemma \ref{l52}, we have a set $V'_{l+1}\subset V_l$, $V'_{l+1}\in\omega$
and $N'_{l+1}\in\rmod$ such that the module
$N'_{l+1}$ $\alpha_{l+1}/2$-tiles $M_j$ if $j\in V'_{l+1}$, where $\alpha_{l+1}\leq \frac{\alpha_l}{2 \dim_K(N_l)}$.
Then by Lemma \ref{l4aug22}, there exists some $m\geq 1$ such that
$N_l$ $\alpha_l$-tiles $(N'_{l+1})^m$. Hence,
$N_l$ $\alpha_l$-tiles $(N'_{l+1})^{m\dim_K(N_l)}$ as well.
Let $N_{l+1}=(N'_{l+1})^{m\dim_K(N_l)}$. By Lemma \ref{l53}, there
exits $n_j>0$ such that
$N_{l+1}$ $\alpha_{l+1}$-tiles $M_j$ provided that $j\in V'_{l+1}$ and
$j>n_j$. So, let 
$$V_{l+1}=V'_{l+1}\cap\{n\,\mid\,n>n_j\}\,.$$
\noindent
Clearly, $V_{l+1}\in\omega$. Then, we satisfy the inductional
hypothesis with $V_1\supseteq V_2 \supseteq \dots\supseteq V_{l+1}$ and the
modules $N_1,N_2,\dots, N_{l+1}$. In order to finish the proof of
our proposition, we need to show that 
for any $k,l\geq 1$ and matrix $A\in\mat_{k\times l}(R)$
$$\lim_{n\to\infty} \rk_{M_n} (A)=\lim_{i\to\infty} \rk_{N_i} (A)\,.$$
\noindent
Therefore, it is enough to prove the following approximation lemma.
\begin{lemma}\label{decappro}
Let $N\subset M$ be finitely generated $R$-modules such that
$\dim_K(N)\geq \dim_K(M) (1-\e)$ for some $0<\e<1$.
Let $k,l\geq 1$ and $A\in\mat_{k\times l}(R)\,.$
Then
$$|\rk_M(A)-\rk_N(A)|\leq 2\e l$$.
\end{lemma}
\proof
By definition,
$$\rk_N(A)=\frac{\rank_{N^l}(A)}{\dim_K(N)},\quad
\rk_M(A)=\frac{\rank_{M^l}(A)}{\dim_K(M)}\,,$$
where $A$ is viewed as a linear map from $M^l$ to $M^k$ (and from $N^l$ to
$N^k$) and $\rank_{M^l}(A):=\dim_K\left(\Ima\mid_{M^l}(A)\right)$. Then,
$$|\rk_M(A)-\rk_N(A)|=\left|\frac{\rank_{M^l}(A)}{\dim_K(M)}-
\frac{\rank_{N^l}(A)}{\dim_K(N)}\right|\leq $$ 
$$\leq \left|\frac{\rank_{M^l}(A)-\rank_{N^l}(A) }{\dim_K(M)}\right|+
\left|\frac{\rank_{N^l}(A)}{\dim_K(M)}-\frac{\rank_{N^l}(A)}{\dim_K(N)}
\right|\leq 2\e l\,$$
\noindent
Hence our lemma and the proposition follows.\qed
\vskip 0.1in
\noindent
Now we turn to the proof of Theorem \ref{tetel2}. We will use the spaces and constants of Proposition \ref{protech}. Set
$n_i=\dim_K(N_i)$, $m_j=\dim_K(M_j)$.
We need three maps for the proof of our theorem, $\rho_\omega:R\to \prod_\omega \Endo_K(M_j)=\mat^\cM_\omega$ is already constructed in Section \ref{conultra}.
\vskip 0.2in
\noindent
\underline{The map $\rho:R\to M_K$:}
\vskip 0.1in
\noindent
Let $V,W$ be finite dimensional $K$-spaces such that $\dim_K(W)=l\dim_K(V)$ and
let $W=W_1\oplus W_2\oplus \dots \oplus W_l$ be a decomposition together with isomorphisms $s_j: V\to W_j$, $1\leq j\leq l$. Then, we have the corresponding
diagonal homomorphism $D:\Endo_K(V)\to \Endo_K(W)$ defined by
$$D(A)(w_1\oplus w_2\oplus\dots\oplus w_l)=$$ $$=\left(s_1(A(s_1^{-1}(w_1)))\oplus s_2(A(s_2^{-1}(w_2)))\oplus\dots\oplus s_l(A(s_l^{-1}(w_l)))\right)\,.$$
Let $N_{i+1}=N_i^{l_i}\oplus V_{i,1}\oplus\dots\oplus V_{i,b_i}$, where the first component is the submodule defined in Proposition \ref{protech} and
$\{V_{i,j}\}^{b_i}_{j=1}$ are arbitrary subspaces satisfying $\dim_K(V_{i,j})=n_i$. Note that the subspaces $V_{i.j}$ can be constructed by the divisibility condition and
that $l_i n_i/n_{i+1}\geq 1-\alpha_i$.
Let $\tau_i:\Endo_K(N_i)\to\Endo_K(N_{i+1})$ be the corresponding diagonal map. So, we have a sequence of injective maps
$$\Endo_K(N_1)\stackrel{\tau_1}{\to}\Endo_K(N_2)\stackrel{\tau_2}{\to}\dots$$
Hence we have a direct limit of matrix algebras
$$T=\varinjlim \Endo_K(N_i)\,.$$
\noindent
Let $\overline{T}$  be the closure of the algebra with respect to the unique rank metric, where $d(A,B)=\Rank(A-B)$.
Then by \cite{Halperin}, $\overline{T}\cong M_K$. Note that this result is originally due to von Neumann. 
Let $\rho_i:R\to \Endo_K(N_i)$ be the map given by the $R$-module structure. Observe that
$$\rk_{N_{i+1}}(\rho_{i+1}(r)-\tau_i\circ\rho_i(r))\leq \alpha_i\,.$$
Hence, the sequence $\{\rho_i(r)\}^\infty_{r=1}$ is Cauchy in $\overline{T}$. Thus, $\lim_{i\to\infty}\rho_i(r)=\rho(r)$ defines a homomorphism
$\rho:R\to\overline{T}\cong M_K\,.$
\vskip 0.2in
\noindent
\underline{The map $\Phi:M_K\to \prod_\omega \Endo_K(M_j)\,$:}
\vskip 0.1in
\noindent
Let $j\in V_i$ and $M_j\cong N_i^{t_{i,j}}\oplus Z_{i,j}$ a decomposition into linear subspaces, where 
the first component $N_i^{t_{i,j}}$ is the submodule obtained by the $\alpha_i$-tiling of $M_j$ by the
module $N_i$ and $Z_{i,j}$ is an arbitrary complementing subspace. So, we have
\begin{equation} \label{egyalfa}
\dim_K(Z_{i,j})\leq \alpha_i\dim_K(M_j)\,.
\end{equation}
\noindent
Using the decomposition above, we define a non-unital injective homomorphism
$\Phi_{i,j}:\Endo_K(N_i)\to\Endo_K(M_j)$ by
$$\Phi_{i,j}(A)=A\oplus A\oplus\dots\oplus A\oplus 0\,,$$
where there are $t_{i,j}$ copies of $A$ in $\Phi_{i,j}(A)$. 
For $j\in\N$, let 
\begin{itemize}
\item $q(j)=j$ if $j\in V_j$.
\item $q(j)=\max\{i\mid\,j\in V_i\}\,,$ if $j\notin V_j$.
\end{itemize}
\noindent
Let $A\in T$ such that $A\in \Endo_K(N_i)$ but $A\notin \tau_{i-1}(\Endo_K(N_{i-1}))$
Then, let
\begin{itemize}
\item $\Psi_j(A)=0$ if $j\notin V_i$.
\item $\Psi_j(A)=\Phi_{q(j),j}(\tau_{q(j)-1}\circ\tau_{q(j)-2}\circ\dots\circ\tau_i(A))$
if $j\in V_i$, that is $q(j)\geq i$.
\end{itemize}
\noindent
Finally,
let
$$\Phi'(A)=[\{\Psi_j(A)\}^\infty_{j=1}]\in\prod_\omega \Endo_K(M_j)\,.$$
\begin{lemma}
$\Phi'$ defines a unital, rank preserving homomorphism from $T$ to $\prod_\omega \Endo_K(M_j)$.
\end{lemma}
\proof
Let $A\in \Endo_K(N_s)$, $A\notin \tau_{s-1}(\Endo_K(N_{s-1}))$ and
$B\in \Endo_K(N_t)$, $B\notin \tau_{t-1}(\Endo_K(N_{t-1}))$.
Then, by Proposition \ref{protech},
$$\{j\in \mid \Psi_j(AB)=\Psi_j(A)\Psi_j(B)\}\in\omega\,$$
and
$$\{j\in \mid \Psi_j(A+B)=\Psi_j(A)+\Psi_j(B)\}\in\omega\,.$$
\noindent
Since the direct limit ring has a unique rank, the only thing remains to be proved is that
$\Phi'$ is unital.
By (\ref{egyalfa}),
$$\rk_{M_j}(\Psi_j(1)-\Id_{\Endo_K(M_j)})\leq \alpha_{q(j)}\,.$$
Since for any $i\geq 1$,
$$\{j\in\N\mid\,q(j)\geq i\}\in\omega$$
\noindent
we have that
$$\lim_{\omega}(\rk_{M_j}(\Phi_j(1)-\Id_{\Endo_K(M_j)})=0\,.$$
Thus, $\Phi'$ is indeed unital. \qed
\vskip 0.1in
\noindent
Now let $\Phi:M_K\to \prod_\omega \Endo_K(M_j)$ be
the closure of $\Phi'$ that is
$$\Phi(\lim_{n\to\infty} x_n)=\lim_{n\to\infty}\Phi'(x_n)\,.$$
\noindent
Since $\Phi'$ is rank preserving, $\Phi$ is a well-defined homomorphism. 
There, we have three homomorpisms.
\begin{itemize}
\item $\rho:R\to M_K$
\item $\rho_\omega:R\to \prod_\omega \Endo_K(M_j)\,.$
\item $\Phi:M_K\to \prod_\omega \Endo_K(M_j)\,.$
\end{itemize}
\vskip 0.1in
The following lemma implies that $\Ima(\phi)$ contains $\Ima(\rho_\omega)$, hence it completes the proof of Theorem \ref{tetel2}.
\begin{lemma}
$\Phi\circ\rho=\rho_\omega\,.$
\end{lemma}
\proof
Let $\kappa_j: R\to \Endo_K(M_j)$ be the homomorphism given by the module structure.
By (\ref{egyalfa}), if $j\in V_i$ and $r\in R$,
$$\rk_{M_j}(\Phi_{i,j}(\rho_i(r))-\kappa_j(r))\leq \alpha_i\,.$$
\noindent
Therefore, for any $i\geq 1$
$$\rk_\omega(\Phi\circ \rho_i(r)-\rho_\omega(r))\leq \alpha_i\,.$$
\noindent
Since, $\lim_{i\to\infty} \rho_i(r)=\rho(r)$ and $\Phi$ is rank preserving map, $\Phi\circ\rho(r)=\rho_\omega(r)\,.$ \qed

\section{String algebras are of amenable representation types}
The significance of string algebras is due to the fact that
the system of f.g. indecomposable modules over such algebras can explicitely be
described (see e.g. \cite{Prest} and the references therein). 
First, let us recall the notion of a string algebra.
Let $Q$ be a finite quiver with vertex set $Q_0$ and arrow set $Q_1$.
Let $KQ$ be the associated path algebra and $I\triangleleft KQ$ be
an ideal, generated by monomials in the path algebras, satisfying the 
following four conditions.
\begin{enumerate}
\item For any vertex $a\in Q_0$, there are at most
two arrows with source $a$, and at most two arrows with target $a$.
\item For any arrow $\alpha\in Q_1$, there exists at most one $\beta\in Q_1$
such that $s(\beta)=t(\alpha)$ and $\beta\alpha\notin I$.
\item For any arrow $\alpha\in Q_1$, there exists at most one $\beta\in Q_1$
such that $t(\beta)=s(\alpha)$ and $\alpha\beta\notin I$.
\item There exists $q\geq 1$ such that any $Q$-path of length at least $q$ is
inside $I$.
\end{enumerate}
\noindent
\begin{example} $KQ=K\langle x,y \rangle, I=\langle x^m,y^n,xy,yx\}$ is a string
algebra. The path-algebra of the $2$-Kronecker quiver is itself a string
algebra.
\end{example}
\noindent
Let $R=KQ/I$ be a string algebra. First we describe the so-called string
modules over $R$, For every arrow $\alpha\in Q_1, \alpha^{-1}$ will denote its
formal inverse such that $s(\alpha)=t(\alpha^{-1})$ $t(\alpha)=s(\alpha^{-1})$.
We denote by $Q^{-1}_1$ the set of all inverse arrows. The elements of $Q_1\cup Q'_1$ are
called {\it letters}. 
A {\it string} of length $n\geq 1$ is a sequence $C_1C_2\dots C_n$
such that
\begin{enumerate}
\item $t(C_{i+1})=s(C_i)$ if $1\leq i \leq n-1$.
\item $C_i\neq C^{-1}_{i+1}$  if $1\leq i \leq n-1$.
\item No substring $C_jC_{j+1},\dots C_k$ or its inverse
$C_k^{-1}\dots C^{-1}_{j+1}C^{-1}_{j}$ lies in the ideal $I$.
\end{enumerate}
\noindent
Additionally, any vertex $a\in Q_0$ is considered to be a string of length
$0$.
Note that for any string $C_1C_2\dots C_n$, $(C_1C_2\dots C_n)^{-1}=
C_n^{-1}C_{n-1}^{-1}\dots C_1^{-1}$ is also a string.
The set of strings is denoted by $W$.
For each string $C_1C_2 \dots C_n=S\in W$ we can associate 
an $R$-module $M(S)$ the following 
way. A $K$-basis for $M(S)$ is $z_0,z_1,\dots,z_n$
For a vertex $a\in Q_0$ and $e_a$ (path of length zero starting and ending
at $a$):
\vskip 0.2\in
\noindent
$e_az_i =
  \begin{cases}
   z_i\,\, \mbox{if $i>0$ and $s(C_i)=a$}\\
   z_0\,\, \mbox{if $i=0$ and $t(C_1)=a$}\\
   0\,\,\,  \mbox{otherwise}
  \end{cases}$
\vskip 0.2\in
\noindent
For an arrow $\alpha\in Q_1$ we define
\vskip 0.2\in
\noindent
$\alpha z_i =
  \begin{cases}
   z_{i-1}\,\, \mbox{if $C_{i-1}=\alpha$}\\
   z_{i+1}\,\, \mbox{if $C_i=\alpha^{-1}$}\\
   0\,\,\,  \mbox{otherwise}
  \end{cases}$
\vskip 0.2\in
\noindent
Observe that $M(S)\cong M(S^{-1})$. On the other hand, if $S_1\neq S_2$,
 $M(S_1)\ncong M(S_2)$
whenever $ S_1\neq S_2^{-1}$. 
\vskip 0.2\in
\noindent
 Now we describe the band modules over the string algebra
$R$.
A string $S=S_{(0)}=C_1C_2\dots C_n$ is called {\it cyclic} if
$S_{(1)}=C_2C_3\dots C_1$,$S_{(2)}=C_3C_4\dots C_2$,
\dots,$S_{(n-1)}=C_nC_1\dots C_{n-1}$ are also strings. Let $S=C_1C_2\dots C_n$
be a cyclic string and $V$ be a finite dimensional vector space over $K$.
Fix an indecomposable transformation $\phi$ on $V$ by identifying
$V$ with $K[x]/f(x)^n$, where $f$ is an irreducible monic polynomial, 
$d=n\,\mbox{deg}(f)$, and $\phi$ is the multiplication by $x$.
The band module $M(S,\phi)$ is defined the following way.
$M(S,\phi)=\oplus^n_{i=1} V_i$, where $V_i\cong V$, that
is for each $v\in V$ and $1\leq i \leq n$ we fix $v_i\in V_i$.
For any $\alpha\in Q_1$ and $i\neq 1,2$:
\vskip 0.2\in
\noindent
$\alpha v_i =
  \begin{cases}
   v_{i-1}\,\, \mbox{if $C_{i}=\alpha$}\\
   v_{i+1}\,\, \mbox{if $C_{i+1}=\alpha^{-1}$}\\
   0\,\,\,  \mbox{otherwise}
  \end{cases}$
\vskip 0.2\in
\noindent
Also,

\noindent
$\alpha v_1 =
  \begin{cases}
   v_n\,\, \mbox{if $C_{1}=\alpha$}\\
   (\phi^{-1}(v))_2\,\, \mbox{if $C_2=\alpha^{-1}$}\\
   0\,\,\,  \mbox{otherwise}
  \end{cases}$
\vskip 0.2\in
\noindent
$\alpha v_2 =
  \begin{cases}
   (\phi(v))_1\,\, \mbox{if $C_{2}=\alpha$}\\
  v_3\,\, \mbox{if $C_3=\alpha^{-1}$}\\
   0\,\,\,  \mbox{otherwise}
  \end{cases}$
\vskip 0.2\in
\noindent
By the Butler-Ringel classification of indecomposable modules \cite{Butler}, $M(S,\phi)\cong M(T,\phi')$ if and only if
$\phi=\phi'$ and $S_{(i)}=T$ for some
$i\geq 1$. Also, a band module and a string module is never isomorphic and any
finitely generated indecomposable module over a string algebra is either
a string module or a band module.
\begin{propo}\label{stringhyper}
Any string algebra $R$ is of hyperfinite type.
\end{propo}

\noindent
The proof will be given in two lemmas.
\begin{lemma}
The family of string modules over a string algebra $R$ is hyperfinite.
\end{lemma}
\proof
By \ref{technical}, it is enough to show that there exists some $m_\e>0$ such that
if $n\geq m_\e$ and $w=C_1C_2\dots C_n$ is a string, then there exist
independent subspaces $W_0,W_1,\dots,W_t\subset M(w)$ such that
\begin{itemize}
\item
$$\dim_K(W_i\oplus \bigoplus_{\alpha_\in Q_1} \alpha W_i)\leq (1+\e) \dim_K(W_i)$$
\noindent
(we only need to check the elements $\alpha_\in Q_1$ not the whole 
path algebra by 
Lemma \ref{linear})
\item 
$$\dim_K(\oplus^t_{i=0} W_i)\geq (1-\e)\dim_K(M(w))\,.$$
\end{itemize}
\noindent
Let $m$ be an integer such that
$\frac{m-2}{m}<1+\e$ and let $m_\e\geq \frac{m}{\e}$. Finally, let
$t\geq 1$ such that $tm\leq n < (t+1)m\,.$

\noindent
Let $W_i$ be the subspace spanned by the vectors 
$\{z_{im},z_{im+1},\dots,z_{(i+1)m-1}\}\,.$
Then,
$$\dim_K(W_i\oplus \bigoplus_{\alpha_\in Q_1} 
\alpha W_i)\leq m+2\leq(1+\e)\dim_K(W_i)\,.$$
\noindent
Also,
$$\dim_K(M(w))-\dim_K(\oplus^t_{i=0} W_i)\leq m\leq \e \dim_K(M(w))\,.\quad\qed$$
\noindent
\begin{lemma}\label{bandlemma}
The family of band modules over a string algebra $R$ is hyperfinite.
\end{lemma}
\proof
The proof will be very similar to the previous one. Fix $\e>0$. It is enough
to prove that there exists $L_\e>0$, such that for any band module
$M(S,\phi)$ there exist independent subspaces 
$\{V_\rho\}_{\rho\in G}$ such that
\begin{itemize}
\item For all $\rho\in G$, $\dim_K(V_\rho)\leq L_\e$.
\item $\dim_K(\sum_{\alpha\in Q_1}\alpha V_\rho +V_\rho)\leq (1+\e)\dim_K(V_\rho)$
\item $\sum_{\rho\in G} \dim_K(V_\rho)\leq (1-\e)\dim_K M(S,\phi)$
\end{itemize}
\noindent
Let $m\in \N$ such that $\frac{m+2}{m}\leq 1+\e$ and
$m_\e>\frac{2m}{\e}$. Let $S=C_1C_2\dots C_n$. First, suppose that $n > m_\e$.
Consider the basis $\{1,x,x^2,\dots,x^{d-1}\}$ for $V$ and let $t\in \N$\
such that $tm+2\leq n < (t+1)m+1\,.$ Then, for any pair
$0\leq  i \leq t-1$ and $0\leq j\leq d-1$ let $W^j_i$ be the subspace
of $M(S,\phi)$ spanned by the set $\{x^j_{im+2},x^j_{im+3},\dots,x^j_{(i+1)m+1}\}$.
Then, it is easy to see that for any pair $i,j$
$$\dim_K(\sum_{\alpha\in Q_1}\alpha W^i_j +W^i_j)\leq (1+\e)\dim_K(W^j_i)\,.$$
\noindent
and
$$\sum_{0\leq  i \leq t-1}\sum_{0\leq j\leq d-1} 
\dim_K(W^j_i)\leq (1-\e)\dim_K M(S,\phi)\,.$$
\vskip 0.2in
\noindent
Now, suppose that for the length of the cyclic string, we have
$n\leq m_\e$ and for the dimension of $V$ we have $d>m_\e$. 
Let $t$ be as above, and
for $0\leq q \leq t-1$ let $Z_q\in M(S,\phi)$ be the subspace spanned by
the set
$$\{\cup_{i=1}^n \cup_{j=qm}^{(q+1)m-1} x^j_i\}\,.$$
\noindent
By the definition of the band module structure,
$$\dim_K(\sum_{\alpha\in Q_1}\alpha Z_q +Z_q)\leq (1+\e)\dim_K(Z_q)\,.$$
\noindent
and
$$\sum_{q=0}^{t-1}\dim_K Z_q\geq (1-\e)\dim_K M(S,\phi)\,.$$
\noindent
Also, for any $q$, $\dim_K(Z_q)\leq mm_\e$. Observe that there are only
finitely many band modules $M(S,\phi)$ for which $n,d\leq m_\e$.
Let $M$ be the maximal $K$-dimension of these modules. So, in order to
finish the proof of the lemma we need to set $L_\e:=\max(mm_\e,M)\,.\,\,$\qed

\section{Benjamini-Schramm convergence implies the convergence 
of string modules}
For the next three sections we fix a string algebra $R$ over the 
finite field $K$. The sole purpose
of the next two sections is to establish a relation between the 
convergence of string modules and
the Benjamini-Schramm convergence of the associated edge-colored graphs. 
This will be the preparation
for the proof of Theorem \ref{stringtheorem}. Let $M(S)$ be the string module corresponding
to the string $S=(C_1C_2\dots C_n)$. We can associate a graph $G_S$ to $S$ 
in a very natural way.
\begin{itemize}
\item The vertex set $V(G_S)$ has $n$ elements $\{x_0,x_1,\dots,x_n\}$.
\item For any $1\leq i \leq n$, we have a directed edge
between $x_{i-1}$ and $x_i$. If $C_i=\alpha\in Q_1$, then
the edge is directed from $x_{i-1}$ to $x_i$ and colored by $\alpha$.
 If $C_i=\alpha^{-1}\in Q_1$, then
the edge is directed from $x_i$ to $x_{i-1}$ and colored by $\alpha$. 
That is, all the
edge-colors are from $Q_1$.
\end{itemize}
\noindent
Note that the graphs $G_S$ and $G_T$ are isomorphic as edge-colored, 
directed graphs if $S=T$
or $S=T^{-1}$. By definition,
if $x$ is a vertex of $G_S$ and $e,f$ are edges both pointing into resp.
 pointing out of $x$, then
the edge-colors of $e$ and $f$ are different.
\begin{defin}
A directed $Q_1$-edge colored graph is called a string graph if all of its
components are in the form of $G_S$ for some string $S$. So, for each 
finitely generated $R$-module
$M$ that is a sum of string modules, we have a unique string graph $G_M$ 
and two $R$-modules
are isomorphic if and only if the corresponding string graphs are isomorphic.
\end{defin}
\noindent
Now we recall the Benjamini-Schramm graph convergence (see \cite{BS} and 
\cite{Elekhyper}). A
rooted string graph is a connected string graph with a distinguished
 vertex $x$ (the root).
The radius of a rooted string graph is the shortest path distance
 between a vertex of $H$ and
the root.
The finite set of rooted string graphs of radius less or equal than 
$r$ is denoted by $U^r$,
that is $U^{r-1}\subset U^r$.
Let $H\in U^r$ and $G$ be a string graph and let $P(H,G)$ denote the set of 
vertices $y$ in $G$
such that the rooted $r$-neighborhood of $y$ is isomorphic to $H$ 
(as rooted, edge-colored, directed graphs). Let $p(H,G):=\frac{|P(H,G)|}{|V(G)|}$.
That is, $p(H,G)$ is the probability that a randomly chosen vertex of $G$
has rooted $r$-neighborhood isomorphic to $H$. The following definition
was originally given for simple graphs \cite{BS}.
\begin{defin}
A sequence of string graphs $\{G_n\}^\infty_{n=1}$ is convergent (in the
sense of Benjamini and Schramm) if for any
$r\geq 1$ and $H\in U^r$, $\lim_{n\to\infty}p(H,G_n)$ exists. \end{defin}
\noindent
Note that for any $H,G$ and $k\geq 1$, $p(H,G)=p(H,G^k)$, where $G^k$ is the disjoint union of $k$ copies of $G$ and it is easy to
see that $p(H,G_1)=p(H,G_2)$ if and only if $G^m_1\cong G^n_2$ for some $m,n\geq 1$.
So, by Proposition \ref{isomodules}, we have the following lemma.
\begin{lemma} Let $M,N$ be finite direct sums of string modules. They represent the
same element in the rank spectrum if and only if 
$p(H,G_M)=p(H,G_N)$ holds for all $r\geq 1$ and $H\in U^r$.
\end{lemma}
\noindent
Now we can state our main proposition in the section.
\begin{propo}\label{egyik}
Let $\{M_n\}^\infty_{n=1}\subset \rmod$ be a sequence of modules such
that for any $n\geq 1$, $M_n$ is the sum of string modules. Suppose that the graphs
$\{G_{M_n}\}^\infty_{n=1}$ converges in the sense of Benjamini and Schramm. Then
$\{M_n\}^\infty_{n=1}$ is a convergent sequence of $R$-modules.
\end{propo}
\noindent
\proof
The proof will be given in a series of lemmas. First we need the definition
of $\e$-isomorphism for string graphs. 
\begin{defin} 
The string graphs $G_1$ and $G_2$ are $\e$-isomorphic, if they contain subgraphs
$J_1\subset G_1$ and $J_2\subset G_2$ such that $J_1$ and $J_2$ are isomorphic string
graphs and $|V(J_1)|\geq (1-\e)|V(G_1)|, |V(J_2)|\geq (1-\e)|V(G_2)|$.
\end{defin}
\begin{lemma}\label{stringisoeps}
For any $\e>0$, there exists $\delta>0$ and $n\geq 1$ such that if the string
graphs $G_1$ and $G_2$ has the same amount of vertices and
$|p(H,G_1)-p(H,G_2)|<\delta$ for any $H\in U^r$, $1\leq r \leq n$, then $G_1$ and $G_2$
are $\e$-isomorphic.
\end{lemma}
\proof (of Lemma \ref{stringisoeps}) Our lemma holds for simple planar graphs with a vertex 
degree bound by the main result of \cite{Newman} (see also Theorem 5. in \cite{Eleklinear}).
Our strategy is to reduce our lemma to the Newman-Sohler result, by enconding our
$Q_1$-edge colored, directed graph with a simple planar graph. 
Let $G$ be a string graph. The encoding simple graph $\hat{G}$ is constructed
as follows. The vertex set of $\hat{G}$ consists of the vertex set of $G$
plus two vertices $a_{xy},b_{xy}$ per each edge of $G$.
If $(x,y)$ is an edge of of $G$, then $(x,a_{xy})$, $(a_{x,y},b_{x,y})$, $(b_{x,y},y)$ will
be edges of $\hat{G}$. That is, we substitute the original edge with a path of length
three. For each element of $\alpha\in Q_1$, we choose a planar graph (the edge-color coding
graph) $T_\alpha$ with a distinguished vertex $t_\alpha$, such that each $T_\alpha$ have the same amount of vertices and
all the vertex degrees of $T_\alpha$ are at least three. For each edge $(x,y)$ of $G$
we stick a copy of $T_\alpha$ to $a_{xy}$ (by identifying $a_{xy}$ and the distinguished
vertex $t_\alpha$) if the edge is directed towards $x$ and its color is $\alpha$.
On the other hand, we stick a copy of $T_\alpha$ to $b_{xy}$ if the edge is
directed towards $y$ and its color is $\alpha$. 
Clearly,
\begin{itemize}
\item We can reconstruct $G$ from $\hat{G}$.
\item If $\{G_n\}^\infty_{n=1}$ is convergent then $\{\hat{G_n}\}^\infty_{n=1}$ is a convergent
sequence of simple planar graphs.
\end{itemize}
\noindent
So, by the Newmark-Sohler Theorem, for any $\e'>0$, there exists
$n\geq 1$ and $\delta>0$ such that if $|p(H,G_1)-p(H,G_2)|<\delta$ holds for any 
$H\in U^r$, $1\leq r \leq n$ and $G_1$ and $G_2$ have the same amount of vertices,
than $\hat{G_1}$ and $\hat{G_2}$ are $\e'$-isomorphic. Also, it is easy to see that for
any $\e>0$ there exists $\e'>0$ such that if 
$\hat{G_1}$ and $\hat{G_2}$ are $\e'$-isomorphic, than $G_1$ and $G_2$ are $\e$-isomorphic.
Hence our lemma follows. \qed
\begin{lemma} \label{1uj}
For any $\e>0$, there exists $n\geq 1$ such that if for the string graphs $G_1$ and $G_2$
$$1-\delta\leq \frac{|V(G_1)|}{|V(G_2)|}\leq 1+\delta$$
and
$|p(H,G_1)-p(H,G_2)|<\delta$ for any $H\in U^r$, $1\leq r \leq n$, then $G_1$ and $G_2$
are $\e$-isomorphic.
\end{lemma}
\proof
First, pick a $\delta$ and $n$ to the constant $\e/2$ as in Lemma \ref{stringisoeps}.
Let $0<\delta'\leq \min(\delta/2, \e/10)$ be a constant such that
for any $H\in U^r$, $1\leq r \leq n$,
$$|p(H,G_2)-p(H,G_3)|<\frac{\delta}{2}$$
provided that
$G_2\subset G_3$ and $|V(G_3)|\leq (1+\delta')|V(G_2)|$. The existence of such $\delta'$ can easily be seen.
Now we show that $\delta'$ and $n$ satisfy the condition of our lemma.
We can suppose that $|V(G_2)|\leq |V(G_1)|$.
By adding strings of one single vertex, we can get $G_2\subset G_3$ such that $|V(G_3)|=|V(G_1)|$.
Then $|V(G_3)|\leq (1+\delta')|V(G_2)|$. Hence for  any $H\in U^r$, $1\leq r \leq n$,
$$|p(H,G_1)-p(H,G_3)|<\delta\,.$$
Therefore, by Lemma \ref{stringisoeps}, $G_1$ and $G_2$ are $\e/2$-isomorphic. Since $\delta'<\e/10$, $G_2$ and $G_1$
are $\e$-isomorphic. \qed
\begin{lemma} \label{2uj}
For any $\e>0$, there exists $\delta>0$ such that if $G_M$ and $G_N$ are $\delta$-isomorphic string graphs, then the
modules $M$ and $N$ are $\e$-isomorphic.
\end{lemma}
\noindent
\proof
We need some notations. If $x\in G_M$, then let $\delta_x$ be the corresponding element in $M$. If $x\stackrel{\alpha}
{\to} y$ is an edge, then $\alpha \delta_x=\delta_y.$ Also, if there is no edge pointing out of $x$ colored by $\alpha$,
then $\alpha \delta_x=0\,.$
By the definition of the string algebra, there exists $q>0$ such that any path of the quiver $Q$ of length at least $q$
is in $I$, where $R=KQ/I$. Therefore, if the $q$-neighborhood of a vertex $x\in G_M$ is in a subgraph $L$, then
\begin{equation}\label{inside}
R\delta_x\subseteq \Span\{\delta_y\mid\, y\in V(L)\}\,.
\end{equation}
\noindent
Suppose that $L_1\subset G_M$, $L_2\subset G_N$ are isomorphic string graphs such that
$|V(L_1)|\geq (1-\delta)|V(G_M)|, |V(L_2)|\geq (1-\delta)|V(G_N)|$. We call $x\in V(L_1)$ an inside point if
its $q$-neighborhood is contained in $L_1$. Obviously, if $\delta$ is small enough, then
$|I_1|\geq (1-\e) |V(G_M)|, |I_2|\geq (1-\e) |V(G_N)|$, where $I_1$ resp. $I_2$ are the sets of inside points
in $L_1$ resp. in $L_2$. Since $L_1$ and $L_2$ are isomorphic, the modules
generated by $\{\delta_x\mid\, x\in I_1\}$ and by $\{\delta_x\mid\, x\in I_2\}$ are also isomorphic. Thus our lemma
follows. \qed
\vskip 0.1in
\noindent
Therefore, from the previous lemmas we have the following corollary.
\begin{corol}
For any $\e>0$, there exists $\delta>0$ and $n\geq 1$ such that
if $M$ and $N$ are sums of string modules and
for any $H\in U^r$, $1\leq r \leq n$,
$$|p(H,G_M)-p(H,G_N)|<\delta$$
\noindent
and also, $1-\delta\leq \frac{\dim_K(M)}{\dim_K(N)}\leq 1+\delta$, then
$M$ and $N$ are $\e$-isomorphic. 
\end{corol}
\noindent
Now, by Lemma \ref{decappro}, our proposition follows. \qed

\section{Convergence 
of string modules implies \\ Benjamini-Schramm convergence}

The goal of this section is to prove
the following converse of Proposition \ref{egyik}.
\begin{propo} \label{masik}
Let $\{M_n\}^\infty_{n=1}\subset \rmod$ be a convergent sequence of modules
such that each module $M_n$ is a sum of string modules.
Then $\{G_{M_n}\}^\infty_{n=1}$ is convergent in the sense of Benjamini and Schramm.
\end{propo}
\noindent
\proof
Let $G$ be a string graph and
$S=C_1C_2\dots C_k$ be a string.
Let $L\cong G_S$ be a subgraph in $G$. We have $V(L)=\{l_0,l_1,\dots,l_k\}$, where
\begin{itemize}
\item The edge $(l_{i-1},l_i)$ is directed towards
$l_{i-1}$ and colored by $C_i$, if $C_i\in Q_1$.
\item The edge $(l_{i-1},l_i)$ is directed towards
$l_i$ and colored by $C_i^{-1}$, if $C_i\in Q^{-1}_1$.
\end{itemize}
\noindent
We call $l_k$ the right resp. $l_0$ the left endvertex of $L$.
By the definition of the strings, if $L_1,L_2$ are subgraphs in $G$ isomorphic
to $G_S$ and their right resp. left endvertices coincide, then $L_1$ and $L_2$.
For the string $S$ and the string graph $G$, let $R(S,G)$ be the set of vertices $x$
in $G$ such that $x$ is the right vertex of a substring $L$ of $G$ isomorphic to $G_S$.
Let
$$r(S,G):=\frac{|R(S,G)|}{|V(G)|}\,.$$
We say that $\{G_n\}^\infty_{n=1}$ is {\bf  stringconvergent}
if for any $S$, $\lim_{n\to\infty} r(S,G_n)$ exists.
The following combinatorial lemma is straightforward to prove.
\begin{lemma}\label{stringconvergent}
If $\{G_n\}^\infty_{n=1}$ is stringconvergent, then it is convergent in the
sense of Benjamini and Schramm.
\end{lemma}
\noindent
Let $N\in\rmod$ be a sum of string modules and $G_N$ be its string graph.
That is, for any $x\in V(G_N)$, we have a base element $\delta_x\in N$ and
if $x\stackrel{\alpha}{\to} y$ is an edge of $G_N$, then
$\alpha(\delta_x)=\delta_y\,$
\begin{lemma}[The String Counting Lemma]
For any string $S=C_1C_2\dots C_k$, there exists a $pp$-pair
$\langle \phi_S,\psi_S\rangle$ such that
$$\dim_K(N(\phi_S))-\dim_K(N(\psi_S))=R(S,G_N)$$
holds for any module $N$ that is a sum of string modules.
\end{lemma}
\noindent
\proof
For $1\leq i \leq k$ let $E_i$ be the equation
\begin{itemize}
\item $C_i n_i- n_{i-1}=0$, if $C_i\in Q_1$.
\item $C_i n_{i-1}- n_i=0$, if $C_i\in Q^{-1}_1$.
\end{itemize}
\noindent Also, let $E_0$ be the equation $n_o=0$.
\noindent
Let
$$N(\phi_S)=\{ n_k\in N\mid\, \exists\, n_0,n_1,\dots,n_{k-1}\subset N\,
\mbox{such that $E_1,E_2,\dots, E_k$ holds}\}\,.$$
$$N(\psi_S)=\{ n_k\in N\mid\, \exists\, n_0,n_1,\dots,n_{k-1}\subset N\,
\mbox{such that $E_0,E_1,\dots, E_k$ holds}\}\,.$$
\noindent
Then clearly, $\phi_S\geq \psi_S$.
Let $\pi:N(\phi_S)\to \Span\{\delta_x:x\in R(S,G_N)\}$ be the natural
restriction map.
\begin{lemma}\label{surjective}
The map $\pi$ is surjective and $\Ker(\pi)=N(\psi_S).$
\end{lemma}
\proof Let $x\in R(S,G_N)$ and $L=[x_0,x_1,\dots,x_k], x_k=x$ be the subgraph
of $G_N$ isomorphic to $G_S$. Then for $n_i=\delta_{x_i}$, $0\leq i \leq k$
the equations $E_1,E_2,\dots, E_k$ hold. Hence $\pi$ is surjective.
If $z\in N(\phi_S)$ and the $\delta_x$-coordinate of $z$ is $\lambda$, then
the $\delta_{x_0}$-coordinate of $z$ is $\lambda$ as well. Hence, if $\pi(z)=0$, then
$z\in N(\psi_S)$. \qed
\vskip 0.1in
\noindent
So we have $$r(S,M_n)=\frac{\dim_K(M_n(\phi_S))}{\dim_K(M_n)}-
\frac{\dim_K(M_n(\psi_S))}{\dim_K(M_n)}\,.$$
By Corollary \ref{ppcorol}, the sequence $\{r(S,M_n)\}^\infty_{n=1}$ 
converges whenever
the modules $\{M_n\}^\infty_{n=1}$
converge. Hence our proposition follows from Lemma \ref{stringconvergent}. \qed

\section{The proof of Theorem \ref{stringtheorem} }
Let $d_R$ be a metric on $\Syl(R)$ defining
the compact topology. Then, for any $\e>0$ there exists
$\delta>0$ and matrices
$A_1,A_2,\dots,A_n\in\mat(R)$ so that if
for any $1\leq i \leq n$, $|\rk_1(A_i)-\rk_2(A_i)|<\delta$,
then $d_R(\rk_1,\rk_2)<\e\,.$
Conversely,
for any set of matrices
$A_1,A_2,\dots,A_n\in\mat(R)$ and $\delta>0$, there exists
$\e>0$ so that if  $d_R(\rk_1,\rk_2)<\e$, then
for all $1\leq i \leq n$, $|\rk_1(A_i)-\rk_2(A_i)|<\delta$.
Particularly, by Lemma \ref{decappro}, for any $\e>0$ there 
exists $\delta>0$ such that if
$M,N$ are $\delta$-isomorphic, then $d_R(\rk_M,\rk_N)<\e.$
First, we prove a weak version of Theorem \ref{stringtheorem}.
\begin{propo}
For any $\e>0$, there exists $\delta>0$ such that
if $M,N\in\rmod$ are sums of string modules,
$1-\delta\leq \frac{\dim_K(M)}{\dim_K(N)}\leq 1+\delta$
and $d_R(\rk_M,\rk_N)<\delta$, then $M$ and $N$ are $\e$-isomorphic.
\end{propo}
\proof
We proceed by contradiction. Suppose that $\{M_n\}^\infty_{n=1}$, $\{N_n\}^\infty_{n=1}$ are sums
of string modules such that $\lim_{n\to\infty}  \frac{\dim_K(M_n)}{\dim_K(N_n)}=1$ and 
none of the pairs $M_n,N_n$ are $\e$-isomorphic. We can also assume that $\{M_n\}^\infty_{n=1}$, $\{N_n\}^\infty_{n=1}$ are
convergent sequences. By Proposition \ref{masik}, 
 $\{G_{M_n}\}^\infty_{n=1}$ and $\{G_{N_n}\}^\infty_{n=1}$ are convergent in the sense of Benjamini and Schramm.
Hence, by Lemma \ref{1uj}, for any $\delta>0$ there exists $n_\delta$ such that
for any $n\geq n_\delta$, $G_{M_n}$ and $G_{N_n}$ are $\delta$-isomorphic. So, by Lemma \ref{2uj},
$M_n$ and $N_n$ are always $\e$-isomorphic, if $n$ is large enough, leading to a contradiction. \qed
\vskip 0.2in
Now suppose that Theorem \ref{stringtheorem} does not hold.
Again, we have two convergent sequence of modules
$\{M_n\}^\infty_{n=1}$ and  $\{N_n\}^\infty_{n=1}$ such that
\begin{itemize}
\item $\lim_{n\to\infty} \frac{\dim_K(M_n)}{\dim_K(N_n)}=1\,.$
\item $\lim_{n\to\infty} \rk_{M_n}= \lim_{n\to\infty} \rk_{N_n}\,.$
\item For each $n\geq 1$,  $M_n,N_n$ are not $\e$-isomorphic.
\end{itemize}
In Lemma \ref{bandlemma} we observed that for any $\kappa>0$, there exists $T_\kappa\geq 0$ such
that if for a band module $B$, $\dim_K(B)\geq T_\kappa$ holds, then there is a module $L_B$ such that
\begin{itemize}
\item $L_B$ is a sum of string modules.
\item $B$ and $L_B$ are $\kappa$-isomorphic.
\end{itemize}
\noindent
 Now, we fix some constants.
Let $\kappa_1>0$ such that if $L_1,L_2$ are sums of string modules and $d_R(L_1,L_2)<\kappa_1$, then
$L_1$ and $L_2$ are $\e/100$-isomorphic. Let $\kappa_2>0$ be a constant such that $\kappa_2<\e/100$ and
if $M,N\in\rmod$ are $\kappa_2$-isomorphic, then $d_R(M,N)\leq \kappa_1/3$.
For each $n\geq 1$, we consider the decompositions,
$$M_n=M^1_n\oplus M_n^2 \quad N_n=N^1_n\oplus N_n^2 \,,$$
where $M^1_n$, $N^1_n$ are sums of band modules of dimension less than $T_{\kappa_2}$ and
$M^2_n,N^2_n$ are sums of string modules and band modules of dimension greater or equal than  $T_{\kappa_2}$.
\begin{lemma} The sequences $\{M^1_n\}^\infty_{n=1}, \{M^2_n\}^\infty_{n=1}, \{N^1_n\}^\infty_{n=1}$ and $\{N^2_n\}^\infty_{n=1}$
are all convergent. Also, the limits \\
$\lim_{n\to\infty}\frac{\dim_K(M^1_n)}{\dim_K(M_n)}, \lim_{n\to\infty}\frac{\dim_K(M^2_n)}{\dim_K(M_n)},
\lim_{n\to\infty}\frac{\dim_K(N^1_n)}{\dim_K(N_n)}$ and $\lim_{n\to\infty}\frac{\dim_K(N^2_n)}{\dim_K(N_n)}$
 exist.
\end{lemma}
\noindent
Le $B_1,B_2,\dots,B_s$ be the set of band modules in $\rmod$ that have dimension less than $T_{\kappa_2}$.
By Corollary \ref{weightcorollary}, for any $1\leq i\leq s$, the limit
$\lim_{n\to\infty} w_{B_i}(M_n)$ exists.
By the definition of the rank (see also Lemma \ref{sulyok}) if $P_1,P_2,\dots,P_l\in\rmod$, then for
any $A\in\mat(R)$,
\begin{equation}\label{teahaz1}
\rk_{\oplus^l_{i=1}P_i}(A)=\sum^l_{i=1} \frac{\dim_K(P_i)}{\dim_K(\oplus^l_{i=1}P_i)} \rk_{P_i}(A)\,.
\end{equation}
\noindent
Recall that
\begin{equation}\label{teahaz2}
M^1_n=\oplus^s_{i=1} B_i^{\frac{w_{B_i}(M_n)\dim_K(M_n)}{\dim_K(B_i)}}\,.
\end{equation}
\noindent
Therefore, $\frac{\dim_K(M^1_n)}{\dim_K(M_n)}=\sum^s_{i=1} w_{B_i}(M_n)$.
Then by (\ref{teahaz2}), $\{M^1_n\}^\infty_{n=1}$ and $\{N^1_n\}^\infty_{n=1}$ are convergent, hence by
(\ref{teahaz1}), $\{M^2_n\}^\infty_{n=1}$ and $ \{N^2_n\}^\infty_{n=1}$ are convergent as well. \qed

\noindent
Also, for large $n$, $B_i^{\frac{w_{B_i}(M_n)\dim_K(M_n)}{\dim_K(B_i)}}$ and $B_i^{\frac{w_{B_i}(N_n)\dim_K(N_n)}{\dim_K(B_i)}}$
are $\e$-isomorphic, hence for large $n$, $M^1_n$ and $N^1_n$ are $\e$-isomorphic.
By our assumptions, we have modules
$L_n\subset M^2_n$, $O_n\subset M^2_n$ that are sums of string modules such that for any $n\geq 1$,
$L_n$ and $M^2_n$ resp. $O_n$ and $N^2_n$ are $\kappa_2$-isomorphic. So, by the definition of $\kappa_2$,
we have that $d_R(L_n,M^2_n)\leq \kappa_1/3$,  $d_R(O_n,N^2_n)\leq \kappa_1/3$. By the previous lemma,
if $n$ is large enough, then $d_R(M^2_n,N^2_n)\leq \kappa_1/3$.
Therefore if $n$ is large enough, then $d_R(L_n,O_n)\leq \kappa_1$. Thus by the definition of $\kappa_1$,
$L_n$ and $O_n$ are $\e/100$-isomorphic. Since $\kappa_2<\e/100$,
we can see that $M^2_n$ and $N^2_n$ are $\e$-isomorphic. Since for large $n$, $M^1_n, N^1_n$ are $\e$-isomorphic,
we get that for large $n$, $M_n$ and $N_n$ are $\e$-isomorphic as well, leading to a contradiction. Hence
our theorem follows.

\section{Parameter testing for modules} \label{property}
A module parameter $p$ is a bounded real function on $\rmod$, where $R$ is a
finite dimensional algebra over a finite field $K$. So, if $M\cong N$ then
$p(M)=p(N)$ (see \cite{Elekpara} for parameters of bounded degree graphs).

\noindent
{\bf Examples:}
\begin{itemize}
\item Let $G(M)$ be the smallest generating system of the module $M$.
Then $g(M):=\frac{|G(M)|}{\dim_K(M)}$ is a module parameter. This parameter
is analogous to the covering number of finite graphs.
\item Let $I(M)$ be the largest system $\{m_1,m_2,\dots, m_{I(M)}\}$ 
of elements in the module $M$ so that $\sum_{i=1}^{|I(M)|} r_im_i=0$ implies that
$r_i=0$ for any $1\leq i \leq |I(M)|$. Then $i(M):=\frac{|I(M)|}{\dim_K(M)}$ 
is a module parameter, analogous to the independence number of a finite graph.
\item Let $Q\in\rmod$. Then the weight function $w_Q(M)$ is a module parameter.
\item Let $Q\in\rmod$. Then $L_Q(M):=\frac{\dim_K(\Hom_R(Q,M))}{\dim_K(M)}$ resp.
$R_Q(M):=\frac{\dim_K(\Hom_R(M,Q))}{\dim_K(M)}$
left and right homomorphism numbers are module parameters analogous to the
left and right homomorphism numbers of finite graphs.
\end{itemize}
\begin{defin}
The module parameter $p:\rmod\to \R$ is stable if it satisfies the following
two conditions:
\begin{enumerate}
\item For any $\e>0$ there exists a $\delta>0$ such that
if $N\subseteq M$ and $\dim_K(N)\geq (1-\delta)\dim_K(M)$ then $|p(M)-p(N)|<\e$.
\item For any $M\in\rmod$ the limit $\lim_{k\to\infty} p(M^k)$ exists.
\end{enumerate}
\end{defin}
\noindent
By Lemma \ref{l4aug22}, for any matrix $A\in\mat(R)$, the matrix parameter $p_A(M)=\rk_M(A)$ is stable.
\begin{propo}
The parameters $g,i,w_Q, L_Q$ and the parameter $R_Q$ provided that $Q$ is
an injective module are stable parameters.
\end{propo}
\proof
Clearly, $w_Q(M^k)=w_Q(M)$,  $L_Q(M^k)=L_Q(M)$, $R_Q(M^k)=R_Q(M)$. Also,
$|G(M\oplus N)|\leq |G(M)|+|G(N)|$ and $|I(M\oplus N)|\geq |I(M)|+|I(N)|$.
So, by Fekete' Theorem on Subadditive and Superadditive Functions
$$\lim_{k\to\infty} \frac{|G(M^k)|}{k\dim_K(M)}\quad\mbox{and}\quad
\lim_{k\to\infty} \frac{|I(M^k)|}{k\dim_K(M)}$$
exist.
Note that similar observations were made by Cohn on projective module
parameters in \cite{Cohnpara}.
It remains to show that the first condition of stability holds for our
parameters.

\vskip 0.1in
\noindent
\underline{$g:\rmod\to \R$\,.}\quad Let $W\subset M$ be a K-linear subspace such
that $N\oplus W=M$. Clearly, if $\{t_1,t_2,\dots, t_s\}$ is a basis of the space $W$ and
$\{n_1,n_2,\dots,n_q\}$ is a generating system for $N$, then $\{t_1,t_2,\dots,t_s,n_1,n_2,\dots,n_q\}$ is
a generating system for $M$.
Therefore
\begin{equation}\label{12d1}
G(M)\leq G(N)+\delta\dim_K(M)
\end{equation}
if $\dim_K(N)\geq (1-\delta)\dim_K(M)\,.$
Now let $\{m_1,m_2,\dots,m_l\}$ be a generating system for $M$ and $\{w_i\}^l_{i=1}\subset W$ be elements
such that $m_i+w_i\in N$, for any $1\leq i \leq l$. Let $n\in N$ and $n=\sum^l_{i=1} r_i m_i\,.$
Then \begin{equation} \label{ize}
n=\sum^l_{i=1} r_i(m_i+w_i)-\sum^l_{i=1} r_iw_i \end{equation}
hence, $\sum^l_{i=1} r_iw_i\in N$. Let $[W]$ be the $R$-module
generated by $W$, so $\dim_K([W])\leq \dim_K(R)\dim_K(W)$.
Let $\{u_1,u_2,\dots,u_s\}$ be a $K$-basis for $[W]\cap N$. Then by (\ref{ize}),
$\{u_1,u_2,\dots,u_s\}\cup \{m_1+w_1, m_2+w_2, \dots, m_l+w_l\}$ is a generating system for $N$.
That is,
\begin{equation} \label{12d2}
G(N)\leq G(M)+\delta \dim_K(M)\dim_K(R)\,.
\end{equation}
\noindent
Now, (\ref{12d1}) and (\ref{12d2}) imply that the first condition of stability holds for the parameter $g$.
\vskip 0.1in
\noindent
\underline{$i:\rmod\to \R$\,.}\quad Clearly, $I(N)\subset I(M)$. Suppose that we have a sequence $\{N_k\subset M_k\}^\infty_{k=1}$ such that
\begin{itemize}
\item $\dim_K(M_k)\to \infty\,.$
\item $\lim_{k\to\infty} \frac{\dim_K(N_k)}{\dim_K(M_k)}=1$, such that $|i(N_k)-i(M_k)|\geq \e.$
\end{itemize}
\noindent
We can suppose that $I(M_k )\geq\e\dim_K(M_k)$.
Let $Z_1,Z_2,\dots,Z_t$ be the finite set of isomorphism classes of principal $R$-modules.
For each $k$, we have $S_k\subset M_k$, $S_k\cong \bigoplus^t_{i=1} Z_i^{q_i}$ such that $\sum^t_{i=1} q_i=I(M_k)$.
By Lemma \ref{l4aug22}, we have $r_i\leq q_i$ and $T_k\subset N_k$, $T_k=\sum^t_{i=1} Z_i^{r_i}$, $\sum^t_{i=1}r_i\leq I(N_k)$
such that
$$\lim_{k\to\infty} \frac{\dim_K(Z_i)(q_i-r_i)}{\dim_K(M)}=0\,.$$
\noindent
Therefore, $\lim_{k\to\infty} |i(N_k)-i(M_k)|=0\,$ in contradiction with our assumption. This implies
that the first condition of stability holds for the parameter $i$.
\vskip 0.1in
\noindent
\underline{$w_Q:\rmod\to \R$\,.}\quad
Again, suppose that we have a sequence  $\{N_k\subset M_k\}^\infty_{k=1}$ such that
\begin{itemize}
\item $\dim_K(M_k)\to \infty\,.$
\item $\lim_{k\to\infty} \frac{\dim_K(N_k)}{\dim_K(M_k)}=1$ and for any $k\geq 1$ 
\begin{equation} \label{12dec3} |w_Q(M_k)-w_Q(N_k)|\geq \e. \end{equation}
\end{itemize}
\noindent
Taking a subsequence, we can assume that $\{M_k\}^\infty_{k=1}$ is convergent. That is, by Lemma \ref{decappro},
for any matrix $A\in\mat(R),$
$$\lim_{k\to\infty} |\rk_{M_k}(A)-\rk_{N_k}(A)|=0\,.$$
\noindent
By (\ref{weight}),
$$w_Q(M_k)=\frac{D_{M_k}(\phi)-D_{M_k}(\psi)}{D_{Q}(\phi)-D_{Q}(\psi)}\quad 
w_Q(N_k)=\frac{D_{N_k}(\phi)-D_{N_k}(\psi)}{D_{Q}(\phi)-D_{Q}(\psi)}$$
\noindent
that is $\lim_{k\to\infty} |w_Q(M_k)-w_Q(N_k)|=0$ in contradiction with (\ref{12dec3}).
So, we established the first stability condition for the parameter $w_Q$.
\vskip 0.1in
\noindent
\underline{$R_Q:\rmod\to \R$\,.}\quad Since $Q$ is injective, any homomorphism $\phi:N\to Q$ extends
to a homomorphism $\phi':M\to Q$. Therefore,
\begin{equation}\label{12dec4}
\dim_K(\Hom(N,Q))\leq \dim_K(\Hom(M,Q))\,.
\end{equation}
\noindent
Let $\pi:\Hom(M,Q)\to\Hom(N,Q)$ be the restriction map.
Clearly, 
$$\dim_K(\Ker(\pi))\leq (\dim_K(M)-\dim_K(N))\dim_K(Q)\,.$$
that is, \begin{equation} \label{12dec5}
\dim_K(\Hom(M,Q))\leq \dim_K(\Hom(N,Q))+ (\dim_K(M)-\dim_K(N))\dim_K(Q)\,.\end{equation}
\noindent
Now, (\ref{12dec4}) and (\ref{12dec5}) imply the first stability condition for the parameter $R_Q$, provided that
$Q$ is injective.
\vskip 0.1in
\noindent
\underline{$L_Q:\rmod\to \R$\,.}\quad
Let us consider the following equations on $M^Q$.
\begin{itemize}
\item For any pair $q_1,q_2\in Q$
$$ m_{q_1}+m_{q_2}-m_{q_1+q_2}=0\,.$$
\item For any pair $r\in R$, $q\in Q$
$$ m_{rq}-rm_q=0\,.$$
\end{itemize}
\noindent
The solution set of the equations above is isomorphic with $\Hom(Q,M)$. So, the
first stability condition for the parameter $L_Q$ immediately follows from the fact that all matrix parameters 
$p_A(M)=\rk_M(A)$ are stable.
\vskip 0.2in
\noindent
The theory of constant-time graph algorithms was developed in the last decade (see e.g. \cite{Gold}). Let us briefly
recall the main idea. Say, we want to estimate the value of a certain graph parameter $p$ for an immensely large
graph $G$ of small vertex degrees. For certain parameters (such as the matching number) we can do the estimate
in constant-time (that is, independently of the size of the graph). A parameter $p$ is called {\bf testable}
for the class of graphs $\cal{G}$, if for any $\e>0$ there exists $\delta>0$ and a family of test-graphs
$F_1, F_2,\dots, F_t$ such that if we learn all the $F_i$-subgraph densities for $G\in\cal{G}$
up to an error of $\delta$, we can
compute $p(G)$ up to an error of $\e$. The point is that using a fixed amount of random samplings of $G$, we can estimate
all the subgraph densities above up to an error of $\delta$ with very high probability, no matter how large our graph
$G$ is. We have an analogous definition for modules.
\begin{defin}
A module parameter $p:\rmod\to \R$ is testable if for any $\e>0$, there exists $\delta>0$, test-matrices
$A_1,A_2,\dots, A_t$ and $n\geq 1$, such that if for a modules $M\in\rmod$, $\dim_K(M)\geq n$, 
and we know all the $\rk_M(A_i)$'s
up to an error of $\delta$, then we can compute $p(M)$ up to an error of $\e$.
\end{defin}
\noindent
The following theorem is motivated by the results of Newman-Sohler \cite{Newman} and
Hassidim-Kelner-Nguyen-Onak \cite{Hassidim}.
\begin{theorem}
If Conjecture \ref{pip} holds for an algebra $R$ (e.g. $R$ is a string algebra), then every
stable parameter $p:\rmod\to \R$ is testable.
\end{theorem}
\noindent
\proof
Let $p$ be a stable parameter for $R$. By our assumptions, we have matrices $A_1,A_2,\dots, A_t$ and $\kappa>0$
such that if
$|\rk_M(A_i)-\rk_N(A_i)|\leq \kappa$ holds for any $1\leq i \leq t$ and
$$1-\kappa<\frac{\dim_K(M)}{\dim_K(N)}\leq 1+\kappa$$
then $|p(M)-p(N)|\leq \e/3\,.$
We can pick $\lambda>0$ in such a way that
\begin{itemize}
\item $\lambda\leq \kappa$
\item if a module $Q$ $\lambda$-tiles $M$ (see Section \ref{hypersection}), then
$|\rk_Q(A_i)-\rk_M(A_i)|\leq \kappa/2$ holds for any $1\leq i \leq t$.
\end{itemize}
\noindent
In the proof of Lemma \ref{l52} we saw that there exists a finite set of modules $\{Q_1, Q_2,\dots, Q_s\}$ such
that for any $M\in\rmod$ at least one of the $Q_j$'s $\lambda$-tiles $M$.
Pick $n$ so large that if $1\leq j \leq s$ and $\dim_K(Q^l_j)\geq n/2$, then
$\lim_{k\to\infty} p(Q^k)-p(Q^l_j)|\leq \e/3\,.$
We also assume that $\dim_K(Q_j)\leq \kappa n$.
Let $p_i=p(Q^n_j)$, then for each $1\leq j \leq s$
\begin{equation}\label{ucso1}
|\lim_{k\to\infty} p(Q^k_j)-p_i|\leq \e/3.
\end{equation}
\noindent
Our algorithm goes as follows. We set $\delta=\kappa/2$. 
If we learn all the values $\rk_M(A_i)$'s up to an error of $\delta$, we can find
at least one $1\leq j\leq s$ such that
\begin{equation}\label{ucso2}
|\rk_{Q_j}(A_i)-\rk_M(A_i)|\leq \kappa
\end{equation}
holds for all $1\leq i\leq t$. These comparisons are all the computations we do after learning
the estimated values $\rk_M(A_i)$. Hence, the computation time does not depend on the dimension of $M$.
Notice however, the role of $n\geq 1$. In the graph parameter case, all the interesting parameters are additive
with respect to disjoint union. Here we only have the stability property. 

\noindent
Note that we might find more than one $Q_j$'s for which (\ref{ucso2}) holds and we cannot be sure that the
$Q_j$ we choose $\lambda$-tiles $M$.  Nevertheless by our assumption, for some $c>0$
$$ 1-\kappa\leq \frac{\dim_K(M)}{\dim_K(Q_j^c)} \leq 1+\kappa\,,$$
hence by the choice of $\kappa$, $|p(M)-p(Q^c_j)|\leq \e/3$, whenever $\dim_K(M)\geq n$.
Therefore by (\ref{ucso1}), $|p_j-p(M)|\leq \e$. So, $p_j$ estimates $p(M)$ up to an error of $\e$.
\qed
\vskip 0.2in
\noindent
{\bf Remark:} One can easily see that the parameters $L_Q$ and $w_Q$ are testable
for any finite dimensional algebra $R$. Indeed, they can be computed from the rank of some matrices.

\vskip 0.2in
\noindent
g.elek@lancaster.ac.uk

\noindent
Lancaster University

\end{document}